\documentclass[11pt]{article}

\usepackage{amsmath, amsfonts,amssymb, amsthm, euscript,makeidx,color,mathrsfs,latexsym}

\newtheorem{assumption}{Assumption}

\newcommand{\ba}{\begin{array}}
\newcommand{\ea}{\end{array}}
\newcommand{\be}{\begin{equation}}
\newcommand{\ee}{\end{equation}}
\newcommand{\bee}{\begin{equation*}}
\newcommand{\eee}{\end{equation*}}
\newcommand{\bea}{\begin{eqnarray}}
\newcommand{\eea}{\end{eqnarray}}
\newcommand{\beaa}{\begin{eqnarray*}}
\newcommand{\eeaa}{\end{eqnarray*}}

\def\dbE{\mathbb{E}}
\def\dbF{\mathbb{F}}

\def\dbL{\mathbb{L}}

\def\dbP{\mathbb{P}}
\def\dbR{\mathbb{R}}

\def\a{\alpha}

\def\d{\delta}
\def\e{\varepsilon}

\def\l{\lambda}

\def\t{\tau}

\def\th{\theta}

\def\D{\Delta}
\def\Th{\Theta}
\def\L{\Lambda}

\def\O{\Omega}
\def\cA{{\cal A}}

\def\cF{{\cal F}}

\def\cL{{\cal L}}

\def\cT{{\cal T}}

\def\cX{{\cal X}}

\def\no{\noindent}

\def\ms{\medskip}
\def\bs{\bigskip}
\def\q{\quad}
\def\qq{\qquad}

\def\pa{\partial}
\def\cd{\cdot}
\def\cds{\cdots}

\def\ol{\overline}
\def\ul{\underline}

\newcommand{\basa}{\begin{assumption}}
\newcommand{\easa}{\end{assumption}}

\newcommand{\bas}{\begin{assum}}
\newcommand{\eas}{\end{assum}}

\def\liminf{\mathop{\underline{\rm lim}}}

\def\pa{\partial}

 \def\cd{\cdot}
\def\cds{\cdots}

\def\dis{\displaystyle}

\def\1{{\bf 1}}

\def\:{\!:\!}
\def\reff#1{{\rm(\ref{#1})}}
\def \proof{{\noindent \bf Proof\quad}}

at 9pt

\begin{document}

\newtheorem{thm}{Theorem}[section]
\newtheorem{lem}[thm]{Lemma}
\newtheorem{cor}[thm]{Corollary}
\newtheorem{prop}[thm]{Proposition}
\newtheorem{rem}[thm]{Remark}
\newtheorem{eg}[thm]{Example}
\newtheorem{defn}[thm]{Definition}
\newtheorem{assum}[thm]{Assumption}

\renewcommand {\theequation}{\arabic{section}.\arabic{equation}}
\def\thesection{\arabic{section}}

\title{\bf A Dynamic Principal Agent Problem with One-sided Commitment}

\author{
Jianfeng Zhang\thanks{\noindent
Department of Mathematics, University of Southern California, Los
Angeles, CA 90089. E-mail: jianfenz@usc.edu. This author is supported in part by NSF grants DMS-1908665 and DMS-2205972.
} ~ and ~{Zimu Zhu}\thanks{\noindent
Department of Statistics and Applied Probability, UC Santa Barbara, Santa Barbara, CA 93106. E-mail: zimuzhu@ucsb.edu. 
}~\thanks{\noindent The authors would like to thank two anonymous referees whose constructive comments have helped to improve the paper greatly.}}

\date{}
\maketitle

\begin{abstract}
In this paper we consider a principal agent problem where the agent is allowed to quit, by incurring a cost. When the current agent quits the job,  the principal will hire a new one, possibly with a different type.  
We characterize the principal's dynamic value function,  which could be discontinuous at the boundary,  as the (unique) minimal solution of an infinite dimensional system of HJB equations,  parametrized by the agent's type. 
This dynamic problem is time consistent in certain sense.

Some interesting findings are worth mentioning. First, self-enforcing contracts are typically suboptimal. The principal would rather let the agent quit and hire a new one. Next, the standard contract for a committed agent may also be suboptimal, due to the presence of different types of agents in our model. The principal may prefer no commitment from the agent, then she can hire a cheaper one from the market at a later time by designing the contract to induce  the current agent to quit. Moreover, due to the cost incurring to the agent, the principal will see only finitely many quittings.

\end{abstract}

\no {\bf Keywords.}   Principal agent problems, contract theory,  one-sided commitment, self-enforcing contracts,  time inconsistency 

\ms

\no{\bf MSC 2020.}  	91B41, 91B43, 93E20, 35K40, 49L25

\eject

\section{Introduction}
The contract theory involves two parties: the principal (she) hires an agent (he) by offering a contract, which provides incentives for the agent to work for her. Such theory has had great impact on organizational economics and corporate finance, as well as many related fields. One typical example, which motivates our work, 
is  employer (principal) versus employee (agent) in labor markets.    Given a contract, the agent's problem is to choose efforts to maximize his expected utility, depending on both the contract and the efforts. Then the principal's problem is to design an optimal contract to maximize her own expected utility, by anticipating that the agent would always choose his corresponding optimal efforts. There is a constraint on the admissible contracts: the agent's optimal utility should satisfy certain {\it individual rationality} (IR) constraint, which can be interpreted as the market value of the agent. The continuous time principal agent problems were first studied by the seminal paper  Holmstr\"{o}m-Milgrom \cite{HM}, followed by Sch\"{a}ttler-Sung \cite{SS} and Sung \cite{S1, S2}. There have been numerous publications on the subject in the past three decades. We refer to the monograph Cvitani\'{c}-Zhang \cite{CZ} and the references therein. In particular, we refer to Sannikov \cite{San}  and Cvitani\'{c}-Possama\"{i}-Touzi \cite{CPT1, CPT2} for the dynamic programming approach, which will be used in this paper.

All the above literature consider the principal agent problem on a fixed time horizon $[0, T]$ (or $[0, \infty)$). There is a crucial time inconsistency issue which {is the underlying reason for a possible renegotiation later} but does not seem to receive much attention in the literature: the optimal contract typically does not remain optimal (or even not admissible anymore) if the two parties reconsider the contract at a later time $t$.\footnote{ We should emphasize that the time inconsistency is not due to the randomness or asymmetric information: even the first best contract in a deterministic setting is typically time inconsistent. Indeed, the optimal contract over period $[t, T]$ should naturally depend on the IR constraint at $t$, which is not involved at all for the problem over $[0, T]$, and thus the optimal contracts over different periods typically do not agree. This exogenously given dynamic IR constraint is the main source of the time inconsistency. In fact, as we see in \cite{CPT1, CPT2, San}, one may choose the (random) IR constraint carefully so as to make the problem time consistent. However, there is no reason to expect this mathematically chosen IR constraint would match the real market values of the agent, and thus the problem in reality is typically time inconsistent.} 
When both parties are fully committed, this issue is irrelevant. However, in many situations, for example in a typical labor contract the agent is allowed to quit before $T$. { Then the principal will need to either renegotiate with the current agent or to hire a new one at the quitting time, with a new contract which is different from the original one exactly due to the time inconsistency.}  Thus, being aware of the agent's possible quitting, the principal would design the contract differently in the beginning.

There have been very serious efforts  on principal agent problems with one-sided commitment. For example,   Grochulski-Zhang \cite{GZ2}, Jeon-Koo-Park \cite{JKK},  Niu-Yang-Zou \cite{NYZ}, Phelan \cite{Phelan},  Ray \cite{Ray}, and Zhang \cite{ZhangYuzhe}  studied models where the agent is not committed. They restricted to the self-enforcing contracts, that is, the agent's continuation utility satisfies the IR constraint at any time.  This is in the sprit of the renegotiation-proof equilibria, see e.g. Farrell-Maskin \cite{FM} and Strulovici \cite{Strulovici}, and under such a contract the agent actually will never quit. We also note that these works mainly consider risk sharing in full information setting, without moral hazard and in particular without the incentive compatibility constraint.   The works Gottlieb-Zhang \cite{GZ} and Karaivanov-Martin \cite{KM} require incentive compatibility, but also restricts the contracts to self-enforcing ones\footnote{It is called non-lapsing constraint in \cite{GZ}, and the so called front-loading constraint is also studied there.}.  The works Ai-Li \cite{AL}, Miao-Zhang \cite{MZ}, and Thomas-Worrall \cite{TW}  are also in the realm of self-enforcing contracts,  but with two-sided limited commitment, namely both the principal and the agent can terminate the contract early. In a different direction,  Krueger-Uhlig \cite{KU} considered a model with multiple principals and agents where the agents are not committed to the contract and can switch to another principal. Their main focus is to determine endogenously the outside option of the agent, which amounts to the IR constraint, through competition of the principals. However, there is no discussion on the principal's behavior after the agent quits, which is the main focus of our paper. The work Ai-Kiku-Li-Tong \cite{AKLT} also endogenizes the outside options, through assortative matching. However, their stopping times for terminating the contract are exogenously given, and before that the contracts are still required to be self-enforcing. Another related work is Hu-Ren-Yang \cite{HRY} where one agent can switch among multiple principals. Their main focus is to find the equilibrium (and mean field equilibrium) contract among the  principals. We also refer to Capponi-Frei \cite{CF}, Hajjej-Hillairet-Mnif \cite{HHM}, He-Tan-Zou \cite{HTZ}, and Lin-Ren-Touzi-Yang \cite{LRTY} where the principal and/or the agent can choose a stopping time to terminate the contract,\footnote{In \cite{CF} the agent chooses a density function for certain stopping time, not the stopping time itself.}  but the principal does not hire a new agent. 
{We should note that, even when restricting to self-enforcing contracts so that there is no desire of renegotiation from the agent's side, the problem may still be time inconsistent.}\footnote{Let's consider an extreme case: the IR constraint  is a real number at $t=0$ but is $-\infty$ at $t>0$. Then the standard optimal contract over $[0, T]$ is self-enforcing. However, it may not be optimal over $[t, T]$ for $t>0$, which has no IR constraint now. We should point out though, in this case the desire of renegotiation is from the principal. When the principal is committed to the contract, as we assume in this paper, then this time inconsistency is irrelevant. } 

We remark that, while disincentivizing the agent from quitting, the self-enforcing contracts are typically sub-optimal for the principal, or say, such a contract is too expensive for the principal. She would rather offer a standard contract, which satisfies the IR constraint only at initial time, and be prepared for a possible quitting of the agent. 
In this paper we consider a continuous time principal agent problem with  moral hazard where the agent is allowed to quit. The main feature of our model is that the principal would hire a new agent when the current agent quits the job, and we shall investigate the time consistency issue from the principal's perspective. One particular application of our model is the dynamic firm-worker relationships, cf. Harris-Holmstr\"{o}m \cite{HH}. There is one risk neutral principal and a family of risk averse agents with different types, parametrized by $\th$. Unlike the adverse selection problems which involve only one agent but with unknown (to the principal) type, here the principal knows the types of all agents but hires only one at any particular time. The principal would pay the agent continuously in time, without lump sum payments. The current agent with type $\th$ would quit the job at time $t$ if he prefers the new opportunity from the market, quantified by the dynamic individual rationality $R^\th_t$. If he quits  at $t$, he will incur a cost $c^\th_t$ and the principal will also suffer a loss $c^P_t$. The principal will  hire a new agent, possibly with a different type $\tilde \th$, and the new agent is also allowed to quit before the expiration date $T$. For simplicity, in this paper we assume $R^\th_t, c^\th_t, c^P_t$ are exogenously given and are deterministic.

Given a contract, the agent can choose his efforts and the quitting time, which is a stopping time. We note that we do not keep track of the agent's status once he quits the current job. This agent's problem is a mixed optimal control and stopping problem, and as standard can be solved through a reflected backward SDE (cf. Zhang \cite{Zhang}). We follow the idea of \cite{CPT1, CPT2, San} to rewrite the agent's continuation (optimal) utility process as a forward process, see also Ma-Yong \cite{MY} for similar idea in the contexts of forward backward SDEs. Then the principal's problem becomes a constrained stochastic control problem, where the agent's continuation utility  becomes the state process, and the controls are the contract and the target effort, with a constraint that the agent's terminal utility is zero (since there is no lump sum payment). The agent's optimal quitting time is an intrinsic hitting time, when the agent's continuation utility process hits the lower barrier $\ul L^\th_t:=R^\th_t - c^\th_t$ (for type $\th$ agent). This enables us to use the HJB equation approach to solve the principal's problem.
 
We now discuss the most crucial component of the problem: the boundary condition of the HJB equation on the lower barrier $\ul L^\th$. When the current type $\th$ agent's continuation utility process hits $\ul L^\th_t$ at $t$, the agent would quit, and then the principal would hire a new agent. First, the principal will choose a different type $\tilde \th$, optimal for the principal's utility. Moreover, the new agent is still allowed to quit, and thus the boundary condition at $\ul L^\th_t$ should be the principal's optimal utility at $R^{\tilde \th}_t$ (minus the cost $c^P_t$), which is exactly the solution we are looking for. So we are considering an infinite dimensional system of HJB equations, parametrized by $\th$, where the interaction is through the boundary condition. The last feature that the boundary condition involves the solution itself actually reflects the time consistent nature of our formulation. However, this adds to the difficulty of the problem, since the boundary condition is not exogenously given.

We establish the wellposedness of this HJB system by recursive approximation. Let $V_n(\th; \cd)$ denote the principal's optimal value function when the current agent is type $\th$ and the principal allows at most $n$ quittings in the remaining time period, that is, after seeing $n$ quittings the principal will offer only self-enforcing contracts for the remaining period. Then the boundary condition for $V_n(\th; \cd)$ at $\ul L^\th_t$ will be $\sup_{\tilde \th} V_{n-1}(\tilde \th; t, R^{\tilde \th}_t) - c^P_t$, which is exogenously given in a recursive way. Under some technical conditions, we show that $V_n$ converges uniformly to the unique solution of the HJB system. Another subtle point is that, at the quitting time $t$, while the agent is indifferent on quitting or staying and thus the agent's optimal utility is continuous at $t$, the principal's utility may have a jump due to the quitting, and thus the principal's value function is typically discontinuous at the boundary  $\ul L^\th_t$. Therefore, the boundary condition is actually an inequality and our value function is characterized as the minimal viscosity solution of the HJB system. 

Some findings implied from our results are interesting from practical perspectives. First, self-enforcing contracts are typically only sub-optimal for the principal.  That is, given that the agent is not committed to the contract, it is not desirable for the principal to offer self-enforcing contracts to disincentivize the agent from quitting, she would rather let the agent quit and rehire a new one when that happens. This observation is more or less obvious, and is strongly confirmed by our results. Next, the standard optimal contract (satisfying the IR constraint only at the initial time) may also be only sub-optimal. That is, in some markets,  the principal may prefer no commitment from the agent so that she can hire a cheaper one from the market at a later time by designing the contract to induce the current agent to quit. The reason is that $R^\th$ (and $c^\th, c^P$) evolves along the time. While it is optimal for the principal to hire an agent $\th$ at time $0$, at some later time $t$ another agent $\tilde \th$ may become much cheaper, and thus the principal would be happy to see agent $\th$ to quit so that she could hire agent $\tilde \th$ to replace him\footnote{Some other principals may still prefer agent $\th$ at time $t$, and that's why $R^\th_t$ could remain high in the market. But we do not model competitions among principals in this paper. We refer again to \cite{AKLT, KU}  for endogenous models for $R^\th_t$.  It will be very interesting to combine our work with those models, and we shall leave this for future research.}.  Moreover, by assuming uniform lower bound of the cost $c^\th$, the agents won't quit too frequently and actually the principal will only see finitely many quittings. This is of course consistent with what we observe in practice. 

We shall mention another direction in the literature on the time inconsistency issue, see e.g.  Balbus-Reffett-Wozny \cite{BRW},  Cetemen-Feng-Urgun \cite{CFU}\footnote{In this work both the principal and the agent have time inconsistent preferences, and they consider certain optimal renegotiation-proof contract.}, Djehiche-Helgesson \cite{DH}, Hern\'{a}ndez-Possama\"{i} \cite{HP},  Li-Mu-Yang \cite{LMY1}, Liu-Mu-Yang \cite{LMY}, Wang-Huang-Liu-Zhang \cite{WHLZ}, Yilmaz \cite{Yilmaz}, and \cite{GZ}, to mention a few. Here the preference of the agent (and/or the principal) is time inconsistent,  for example due to non-exponential discounting. We refer to the survey paper Strotz \cite{Str} for time inconsistent stochastic control problems, in particular, one may consider pre-committed, naive, and sophisticated agents. The principal knows what type of agent she is hiring, and would find the optimal contract (or possibly sub-optimal one when the principal's preference is also time inconsistent) by anticipating that the agent would always choose his optimal or sub-optimal efforts in a corresponding way.  We emphasize again that the time inconsistency in these works is due to the individual player's preference. This has completely different nature from the time inconsistency we investigate in this paper due to the interaction between the principal and the agents. 
 
Finally we clarify the time consistency we obtain in this paper. First, our agent's problem is always time consistent for the period he stays with the job, and we do not keep track of the agent once he leaves the current principal. Our principal's problem is time consistent at the agents' optimal quitting times, see Theorem \ref{thm-DPP}. That is, the principal would never regret when an agent quits and thus she could offer a new contract to another agent. We should note though, if the principal reconsiders the contract within the time period an agent is staying, she may find it not optimal. However, since the principal is required to commit to the contract, she is not allowed to fire the agent and sign a better contract during those time periods, and thus such time inconsistency for the principal is irrelevant from practical considerations.

The rest of the paper is organized as follows. In Section \ref{sect-0quit} we introduce our model and consider self-enforcing contracts. In Section \ref{sect-1quit} we consider the model where only one quitting of the agent is allowed, and the principal's dynamic value function in this model is analyzed in Section \ref{sect-u1}. In Section \ref{sect-nquit} we extend these results  to the case with at most $n$ quittings, and in Section \ref{sect-general} we remove the restriction on the number of quittings. In Section \ref{sect-standard0} we reinvestigate the problem when the benchmark model is the standard principal agent problem with full commitment. Finally in Section \ref{sect-summary} we summarize our findings, and in Appendix we complete some technical proofs.

\section{Our model and self-enforcing contracts}
\label{sect-0quit}
\setcounter{equation}{0}
In this section we introduce our model on the fixed time horizon $[0, T]$. In particular, we shall consider self-enforcing contracts so that the agent has no incentive to quit. 

Let $(\O, \dbF, \dbP)$ be a filtered probability space on $[0, T]$, $B$ a standard one dimensional Brownian motion, and $\dbF=\dbF^B$. There are two players: a principal (she) and an agent (he). The agent has a parameter $\th\in \Th$, which can be viewed as the type of the agent. In this section, $\th$ is fixed, while in the rest of the paper, we will consider different types of agents.

The principal's control, namely the contract, is a continuous time payment  $\eta = \{\eta_t\}_{0\le t\le T}$. The agent's control, namely his action or effort, is a process $\a = \{\a_t\}_{0\le t\le T}$.  In particular, there are no terminal payments, which are not convenient to model when an agent quits. The admissible controls satisfy the following  conditions, required for technical convenience.

\begin{defn}
\label{defn-cA}
(i) Let $\cA^P$ denote the set of $\dbF$-progressively measurable processes $\eta$ on $[0, T]$ such that $-C(\eta) \le \eta  \le C_0 $. That is, $\eta$ is bounded, with a uniform upper bound $C_0>0$, but the lower bound $-C(\eta)$ may depend on $\eta$. 

(ii) Let $\cA^A$ denote the set of $\dbF$-progressively measurable processes $\a$ on $[0, T]$ such that
\bea
\label{BMO}
\dbE^\dbP_t\Big[ \int_t^T |\a_s|^2 ds \Big] \le C(\a),\q 0\le t\le T,\q \dbP\mbox{-a.s.}
\eea
where the constant $C(\a)>0$ may depend on $\a$. That is, $\int_0^t \a_s dB_s$ is a BMO martingale.
\end{defn}

 We consider moral hazard models, and as in the standard literature we use the weak formulation as follows, see e.g. Cvitani\'{c}-Wan-Zhang \cite{CWZ}. The state process is:
\bea
\label{X}
X_t \equiv B_t.
\eea
The agent controls the law of $X$, denoted as $\dbP^\a$, through the Girsanov theorem:
\bea
\label{Pa}
B^\a_t := B_t - \int_0^t \a_s ds,\q M^\a_t:= \exp\Big(\int_0^t \a_s dB_s - {1\over 2} \int_0^t |\a_s|^2 ds\Big),\q d\dbP^\a := M^\a_T d\dbP.
\eea
Note that under \reff{BMO}  $M^\a$ is a true martingale under $\dbP$ and hence $\dbP^\a$ is another probability measure. Moreover, there exists $\e=\e(\a)>0$ such that (cf. \cite[Section 7.2]{Zhang})
\bea
\label{BMO2}
\dbE^\dbP\Big[ \exp\big(\e \int_0^T |\a_s|^2 ds\big) + |M^\a_T|^{1+\e}\Big] <\infty.
\eea
We shall also introduce the following notations:
\bea
\label{Ea}
\dbE := \dbE^\dbP,\q \dbE^\a := \dbE^{\dbP^\a},\q \dbE^\a_t:= \dbE^{\a}[~ \cd~ |\cF_t].
\eea

\subsection{The agent's problem}
Given $\eta\in \cA^P$, the agent's  optimal expected utility at time $0$ is:
\bea
\label{JA}
\dis V^A_0(\th; 0, \eta)  := \sup_{\a\in \cA^A} \dbE^{\a}\Big[\int_0^T  f(\th; s,\eta_s) ds-\int_0^T  h(\th; s,\a_s) ds\Big].
\eea

\no Here the agent's running utility function $f$ and running cost function $h$ satisfy Assumption \ref{assum-uAh}  below. 
Moreover, the subscript $_0$ in $V^A_0$ indicates that the agent is allowed to quit $0$ time, namely  not allowed to quit.

\begin{assum}
\label{assum-uAh}
There exists a constant  $\L_0 \ge 1$ satisfying the following.

(i) $f$ is twice differentiable in $\eta$ with:
\bea
\label{uA}
\left.\ba{c}
\dis  f\le -{1\over \L_0}<0,\q \pa_\eta f \ge {1\over \L_0}> 0,\q \pa_{\eta\eta} f \le -{1\over  \L_0}<0,\q \mbox{for all}~ \eta \le C_0;\ms\\
\dis\mbox{and}\q f(\th; s, C_0) \ge -\L_0.
\ea\right.
\eea

(ii) $h$ is  twice differentiable in $\a$ with:
\bea
\label{h}
h\ge 0~\mbox{with}~ h(\th; s,0)= 0,\q {1\over \L_0} \le \pa_{\a\a} h\le \L_0.
\eea

(iii) $f$, $h$, $\pa_\a h$ are uniformly  continuous in $s$, with a common modulus of continuity function $\rho$, and we assume without loss of generality that $\rho(\d)\ge \d$.
\end{assum}
\no { Basically we require the utility function $f$ to be increasing and concave in the payment $\eta$ and the cost function $h$ to be convex in the effort $\a$, which are standard, and the rest of the assumption is for technical convenience.} Clearly, in this section we may allow $\L_0$ and $\rho$ to depend on $\th$. However, for the purpose of later sections, we require them to be common for all $\th$. Throughout the paper, we shall denote by $C$ a generic constant which may depend on $T, C_0$ and $\L_0$, but not on $\rho$ or $\th$. One typical example satisfying Assumption \ref{assum-uAh} is: for some appropriate $\l_1(\th), \l_2(\th)>0$, 
\bea
\label{quadratic}
f(\th; s, \eta) = -\l_1(\th) e^{-\l_1(\th) \eta},\q h(\th; s, \a) = {1\over 2\l_2(\th)} |\a|^2.
\eea

Under \reff{uA} and \reff{h}, it is clear that $f(\th; s, \cd)$ and $\pa_\a h(\th; s,\cd)$ have inverse functions, denoted as $I^f(\th; s,\cd)$ and $I^h(\th; s, \cd)$, respectively. That is, 
\bea
\label{I}
f\big(\th; s, I^f(\th; s, y)\big)=y, ~ y\le f(\th; s, C_0);\q~ \pa_\a h\big(\th; s, I^h(\th; s, z)\big) = z,~ z\in \dbR.
\eea
Moreover, we introduce the following functions:
\bea
\label{H}
\dis H(\th; s, z) :=\sup_{\alpha\in \dbR}\big[\alpha z-h(\th; s,\a)\big],\q  \tilde H(\th; s, z) := -h(s, \th, I^h(\th; s, z)).
\eea
Then  $I^h$ is the maximum argument of $H$: $H(\th; s, z) = zI^h(\th; s, z) -h(s, \th, I^h(\th; s, z))$,
and it is straightforward to verify the following: for some constant $C(\L_0)>0$,
\bea
\label{HIhproperty}
\left.\ba{c}
\dis 0<\pa_y I^f \le \L_0,\q  |\eta|^2 \le  C(\L_0)\Big[1 -f(\th; s, \eta)\Big],\q \mbox{for all}\q \eta \le C_0;\\
\dis I^h(\th; s, 0) = 0,\q 0< {1\over \L_0} \le \pa_z I^h(\th; s, z) \le \L_0 <\infty;\\
\dis I^h~\mbox{is uniformly continuous in $s$ with modulus of continuity function $\L_0\rho$};\ms\\
\dis \tilde H(s, \th, z) \le  \tilde H(s, \th, 0) = 0= H(s, \th, 0) \le H(s, \th, z);\ms\\
\dis |\pa_z H(\th; s, z)|,  |\pa_z\tilde H(\th; s, z)| \le \L_0|z|,\q  |z|^2 \le -C(\L_0)\tilde H(\th; s, z).
\ea\right.
\eea

{
The solution to the agent's problem  is standard, see e.g. \cite{CZ}, which provides both the agent's optimal utility and his optimal action through the following BSDE: given $\eta\in \cA^P$,
\bea
\label{BSDE}
 Y^{\th, \eta}_t =  \int_t^T \big[ f(\th; s,\eta) + H(\th; s,Z^{\th,\eta}_s)\big] ds - \int_t^T Z^{\th,\eta}_s dB_s,\q t\in [0, T].
 \eea

\begin{prop}
\label{prop-agent}
Let Assumption \ref{assum-uAh} hold and $\eta\in \cA^P$.  Then the agent has optimal utility $V^A_0(\th; 0, \eta) = Y^{\th,\eta}_0$, and has the unique optimal control  $\a^*= I^h(\th; \cd,Z^{\th,\eta})\in \cA^A$.
 \end{prop}
 
 Following \cite{MY, San}, we shall rewrite the BSDE \reff{BSDE} in the forward form,  which will be crucial in the next section. For any $t\in [0, T)$, $x\in \dbR$, and $(\eta, Z)\in \cA^P_t\times \cA^A_t$, denote
\bea
\label{Yforward}
\left.\ba{c}
\dis \cX^{\th; t, x, \eta, Z}_s  := x - \int_t^s \big[f(\th; r,\eta_r) + H(\th; r,Z_r)\big] dr + \int_t^s Z_r dB_r \\
\dis \qq\qq =x  - \int_t^s \big[f(\th; r,\eta_r) +\tilde{H}(\th; r,Z_r)\big] dr + \int_t^s Z_r dB_r^{\th, Z},~ t\le s\le T,
\ea\right.
\eea
and set, recalling $Y^{\th,\eta}_T=0$,  
\bea
\label{cAy}
\cA^\th_{t, x} := \Big\{(\eta, Z)\in \cA^P_t\times \cA^A_t: \cX^{\th; t, x, \eta, Z}_T=0\Big\}.
\eea
}

\begin{rem}
\label{rem-forward}
(i) We emphasize that $\cX^{\th; t, x, \eta, Z}$ is the forward form of $Y^{\th, \eta}$, which stands for the continuation utility of the agent, rather than the state process $X$ in \reff{X}. In particular, $x$ is the agent's optimal utility at initial time $t$ for the contract $\eta$.

 (ii) In \cite{CPT1, CPT2}, the principal has terminal payment $\xi$, then $\xi$ has one to one correspondence with the corresponding $(x, Z)$, and thus one may view $(x, Z)$ as the principal's controls. Alternatively (although less natural in this case), we may view $(\xi, x, Z)$ altogether as the principal's control, but then require a constraint $\xi = \cX^{t,x,Z}_T$ for some appropriate process $\cX^{t,x,Z}$ in the spirit of \reff{Yforward}. In our situation, $\eta$ and $(x, Z)$ also have one to one correspondence, however, it is not easy to express $\eta$ directly in terms of $(x, Z)$. So we take the alternative way by viewing the triple $(\eta, x, Z)$ as the princilal's controls, but require the constraint $ \cX^{\th; t, x, \eta, Z}_T=0$ in \reff{cAy}.
 \end{rem}

  Recall \reff{Pa} and \reff{Ea}, we denote
 \bea
 \label{EthZ}
 B^{\th,Z} := B^{I^h(\th;\cd, Z)},\q \dbP^{\th,Z} := \dbP^{I^h(\th;\cd, Z)},\q \dbE^{\th,Z} := \dbE^{I^h(\th;\cd, Z)},\q \dbE_t^{\th,Z} := \dbE_t^{I^h(\th;\cd, Z)}.
 \eea
Moreover, recall Definition \ref{defn-cA} (i) and introduce
  \bea
 \label{Ltht}
\ol L^\th_t:= \int_t^T f(\th; s,C_0) ds,\q 0\le t\le T.
 \eea 
 We note that, when $\eta\equiv C_0$, obviously $(\ol L^\th_t, 0)$ satisfies BSDE \reff{BSDE}, so $(Y^{\th, C_0}_t, Z^{\th,C_0}_t) = (\ol L^\th_t, 0)$.  Since $\eta\le C_0$,  by applying the comparison principle for BSDE \reff{BSDE}, we have 
 \bea
 \label{olLmax}
 Y^{\th, \eta}_t \le \ol L^\th_t ~\mbox{for all}~ \eta \in \cA^P,\q\mbox{and thus},\q  \cA^\th_{t,x}\neq \emptyset~\mbox{if and only if} ~x\le \ol L^\th_t.
 \eea
 The following technical estimate will be used frequently in the paper. We postpone its proof to Appendix. 
 
 \begin{lem}
\label{lem-YZetaest}
 Let  Assumption \ref{assum-uAh}  hold. For any $t<T$ and $\eta \in \cA^P$,   we have
\bea
 \label{YZetaest}
 \left.\ba{c}
\dis \dbE^{ \th,Z^{\th, \eta}}_t\Big[\int_t^T \big[|Z^{\th,\eta}_s|^2 + C_0-\eta_s\big] ds\Big] \le C\big[\ol L^\th_t - Y^{\th,\eta}_t\big],\\ 
\dis \dbE^{ \th,Z^{\th,\eta}}_t\Big[\int_t^T |\eta_s|^2 ds\Big]\le C\big[T-t-Y^{\th,\eta}_t\big].
\ea\right.
 \eea
 Consequently, for any $t<T$, $x\le \ol L^\th_t$, and $(\eta, Z)\in \cA^\th_{t,x}$, we have
 \bea
 \label{Zetaest1}
 \dbE^{ \th,Z}\Big[\int_t^T \big[|Z_s|^2 + C_0-\eta_s\big] ds\Big] \le C[\ol L^\th_t - x],\q \dbE^{ \th,Z}\Big[\int_t^T |\eta_s|^2 ds\Big]\le C[T-t-x].
 \eea
 \end{lem}

 \subsection{The principal's problem}
 We consider a risk neutral principal. Given $\eta\in \cA^P$, the principal's expected utility is: 
 \bea
\label{JP}
J_P(\th; 0, \eta)  := \dbE^{\a^*}\Big[X_T - \int_0^T \eta_sds\Big] 
=  \dbE^{\th, Z^{\th,\eta}}\Big[  \int_0^T  \big[ I^h(\th; s,Z^{\th,\eta}_s)-\eta_s\big] ds\Big].
\eea

\begin{rem}
\label{rem-generalization}
(i) In \reff{Pa} we assume the drift of the state process is $\a$. This is a simplification without loss of generality. Indeed, if we consider a more general drift $b(\th; s, \a_s)$, we may simply set a new control $\tilde \a_s := b(\th; s, \a_s)$, and one can easily see that the agent's problem \reff{JA} is equivalent to the problem by replacing $h$ with $\dis \tilde h(\th; s, \tilde \a_s) := \inf_{\a_s: b(\th; s, \a_s) = \tilde \a_s} h(\th; s, \a_s)$. 

(ii) As in the standard literature for principal agent problems with one-sided commitment, in \reff{JP} we consider a risk neutral principal. We note that the slightly more general case with a discount factor $r$ can be easily transformed into our model. Indeed,  consider
\beaa
J_P(\th; 0, \eta)  := \dbE^{\a^*}\Big[e^{-r T} X_T - \!\!\!\int_0^T\!\!\! e^{-r s} \eta_sds\Big] 
=  e^{-r T} \dbE^{{\th,Z^{\th,\eta}}}\Big[  \int_0^T \!\!\! \big[ I^h(\th; s,Z^{\th,\eta}_s)- e^{r(T- s)}\eta_s\big] ds\Big].
\eeaa
Denote 
\beaa
\tilde \eta_s := e^{r(T- s)}\eta_s,~  \tilde Z^{\th; \tilde \eta} := Z^{\th, \tilde\eta e^{r(\cd-T)}},~ \tilde J_P(\th; 0, \tilde \eta) := \dbE^{\th,\tilde Z^{\th, \tilde\eta}}\Big[  \int_0^T [ I^h(\th; s, \tilde Z^{\th,\tilde\eta}_s)- \tilde \eta_s] ds\Big].
\eeaa
Then the principal's problem is equivalent to optimizing $\tilde J_P$, which has no discount factor. 

(iii) The case with a risk averse principal is technically harder. Moreover, in \reff{JA} we assume the agent's utility $f$ and cost $h$ are separable. The general non-separable case in the form $f(\th; s, \eta, \a)$ is also technically harder. We do not pursue such generality in this paper.
\end{rem}

Introduce the agent's individual rationality $R^\th_t$ at time $t$, which can be interpreted as the agent's market value if he works for a principal over period $[t, T]$.  In this paper we assume $R^\th$ is deterministic and is given exogenously, and we leave more general cases to future research. In light of \reff{olLmax}, we assume the following, with the same $\L_0, \rho$ for simplicity.

\begin{assum}
\label{assum-R}
(i) $-\L_0 \le R^\th_s \le \ol L^\th_s$, $0\le s\le T$, and $R^\th_T =0$.

(ii) $R^\th$ is uniformly Lipschitz continuous in $s$, with modulus of continuity function $\rho$. 
\end{assum}

In the standard literature, the principal's problem is as follows:
\bea
\label{VP0standard}
V^{S,P}_0(\th; 0)  := \sup_{\eta\in \cA^P} J_P(\th; 0, \eta),\q\mbox{subject to}\q V^A_0(\th; 0, \eta) \ge R^\th_0.
\eea
Here $^S$ stands for standard. We shall provide more discussions on $V^{S,P}_0(\th; 0)$ in Section \ref{sect-standard} below. In this paper we mainly focus on self-enforcing contracts: recalling Proposition \ref{prop-agent},
\bea
\label{VP0}
V^P_0(\th; 0)  := \sup_{\eta\in \cA^P} J_P(\th; 0, \eta),\q\mbox{subject to}\q Y^{\th,\eta}_t \ge R^\th_t, ~ 0\le t\le T, \mbox{a.s.}
\eea
By \reff{olLmax} it is clear that $\eta \equiv C_0$ satisfies the constraints and thus
\bea
\label{VP0-C0}
V^{S,P}_0(\th; 0) \ge V^P_0(\th; 0) \ge J_P(\th; 0, C_0) = - C_0(T-t).
\eea

\subsection{The principal's dynamic value}
We next extend the above principal agent problem to arbitrary time interval $[t, T]$. Denote by $B^t_s:= B_s- B_t$ the shifted  Brownian motion on $[t, T]$,  $\dbF^t:= \dbF^{B^t}$, and $\cT_t$  the set of $\dbF^t$-stopping times on $[t, T]$. We shall extend all the notations in the previous subsection to the shifted space in the obvious manner. In particular, $\cA^P_t, \cA^A_t$ will denote the set of $\dbF^t$-progressively measurable processes $\eta, \a$, but we shall still denote by $(Y^{\th, \eta}, Z^{\th, \eta})$ the solution to the BSDE \reff{BSDE} on $[t, T]$. We then define the principal's dynamic value over self-enforcing contracts as follows: 
\bea
\label{VPt}
\left.\ba{c}
\dis V^P_0(\th; t)  := \sup_{\eta\in \cA^P_t} J_P(\th; t, \eta),\q\mbox{subject to}\q Y^{\th, \eta}_s \ge R^\th_s, ~ t\le s\le T, \mbox{a.s.}\\
\dis \mbox{where}\q J_P(\th; t,  \eta):= \dbE^{\th, Z^{\th;\eta}}\Big[  \int_t^T \big[I^h(\th; s,Z^{\th; \eta}_s)-\eta_s\big] ds\Big].
\ea\right.
\eea

{ The following result states that the principal's optimal utility over self-enforcing contracts is uniformly continuous with respect to the contracting period, more precisely the initial time of the contract. This is interesting in its own right, and will be important technically in the next section. We postpone its proof to Appendix.}
\begin{prop}
\label{prop-VP0}
Under Assumptions \ref{assum-uAh} and \ref{assum-R},  $V^P_0(\th; \cd)$ is uniformly continuous in $t$:
\bea
\label{VP0-reg}
|V^P_0(\th; s) - V^P_0(\th; t)|\le C\sqrt{\rho(|s-t|)},\q \forall s, t\in [0, T].
\eea
 \end{prop}

\section{The case with at most one quitting}
\label{sect-1quit}
\setcounter{equation}{0}
We now assume the principal offers a standard contract which meets the individual reservation only at initial time. The agent might quit before $T$, at a stopping time up to his choice, if he finds a better opportunity in the market. Then the principal will hire another agent with a new contract. In this section, we consider the case that, once the current agent quits the job, the principal will offer only self-enforcing new contract. That is, the principal will tolerate only one quitting.  We assume there is a pool of agents in the market with parameter set $\Th$.  For each agent $\th \in \Th$, if he quits at time $t$, he will incur a cost $c^\th_t$, and the principal will bear  a cost $c^P_t$, which (for simplicity) does not depend on $\th$.  Recall \reff{olLmax} and denote further that
\bea
\label{ulLtht}
\ul L^\th_t := R^\th_t - c^\th_t,\qq D_\th := \{(t,x): 0\le t<T, \ul L^\th_t < x \le \ol L^\th_t\},\q T_\d:= T-\d.
\eea
In this paper we assume $c^\th, c^P$ are deterministic and are given exogenously, and we leave more general cases to future research. The following technical conditions will be in force.  

\begin{assum}
\label{assum-RCth}
$R^\th$ satisfies Assumption \ref{assum-R} and, for the same $\L_0, \rho$ (for simplicity), 

(i)  $c^\th_t \ge {1\over \L_0}>0$,  $c^P_t \ge 0$, and consequently, $\ol L^\th_t - \ul L^\th_t \ge {1\over \L_0}$, $0\le t\le T$;

(ii)  $c^\th, c^P$ are uniformly continuous  in $t$ with modulus of continuity function $\rho$;

(iii) $c^\th, c^P$ are  bounded by $\L_0$, and $|\ul L^\th_{t} - \ul L^\th_s |\le \L_1 |t-s|$ for some $\L_1>0$.
\end{assum}
\no From now on, we allow the generic constant $C>0$ to depend on $\L_1$ as well.

\subsection{The agent's problem}
 Consider the problem on time interval $[t, T]$ with an agent $\th$ at initial time $t$ and a contract $\eta \in \cA^P_t$.  The agent can choose a stopping time $\t \in \cT_t$ to quit the job.  In that case, he would receive the payment $\eta$ from the principal over $[t, \t]$, and would get a new job with expected utility $R^{\th}_\t$ for the remaining period $[\t, T]$, and in the mean time he would incur the cost $c^\th_\t$. He may also choose not to quit, and in that case we set $\t = T$. Therefore, given the contract $\eta\in \cA^P_t$, the agent's optimal utility is: 
\bea
\label{VA1}
V^A_1(\th; t, \eta) := \sup_{\t\in \cT_t} \sup_{\a\in \cA^A_t} \dbE^\a\Big[\int_t^\t \big[f(\th; s,\eta_s)-h(\th; s,\a_s)\big]ds + \ul L^\th_\t \1_{\{\t<T\}}\Big].
\eea
Here the subscript $_1$ in $V^A_1$ indicates that the agent can quit once. This problem is solved through a reflected BSDE (cf.  \cite[Chapter 6]{Zhang}): 
\bea
\label{RBSDE}
\left.\ba{c}
\dis Y_s^{\th; t, \eta}=\int_s^T \big[f(\th; r,\eta_r) + H(\th; r,Z^{\th;t,\eta}_r)\big] dr - \int_s^T Z^{\th;t,\eta}_r dB_r+K^{\th;t,\eta}_T-K^{\th;t,\eta}_s,\\
\dis Y^{\th;t,\eta}_s \geq \ul L^{\theta}_s,\qq  \int_t^T [Y^{\th;t,\eta}_s-\ul L^{\theta}_s] ~dK^{\th;t,\eta}_s=0.
\ea\right.
\eea
Here the solution is a triplet of processes $(Y^{\th;t, \eta}, Z^{\th;t, \eta}, K^{\th; t, \eta})$, and $K^{\th; t, \eta}$ is a non-decreasing continuous process with $K^{\th; t, \eta}_t = 0$. In particular, we use the subscript $^t$ to denote the solution to RBSDE \reff{RBSDE}, while using $(Y^{\th, \eta}, Z^{\th, \eta})$ for BSDE \reff{BSDE} (even if it is on $[t, T]$). We introduce the first time that $Y^{\th;t,\eta}$ hits the barrier:
\bea
\label{taueta}
\t^{\th; t,\eta}:= \inf\big\{s\ge t:  Y^{\th; t, \eta}_s = \ul L^{\th}_s\big\} \wedge T.
\eea
We then have the following result which solves the agent's problem in this case.
\begin{prop}
\label{prop-agent1}
Let Assumptions \ref{assum-uAh} and \ref{assum-RCth} hold true. Then, for any $\th\in \Th$ and $t<T$, 

\ms
(i) $V^A_1(\th; t, \eta) = Y^{\th; t, \eta}_t$, the  $\t^{\th; t,\eta}\in \cT_t$ in \reff{taueta} is an optimal quitting time and $\a^{\th; t,\eta} := I^h(\th; s,Z^{\th;t,\eta})\1_{[t,\tau^{\th;t,\eta}]}\in \cA^A_t$ is the agent's unique optimal effort (unique up to $\t^{\th; t,\eta}$).

\ms
(ii) $\{Y^{\th; t, \eta}_t: \eta\in \cA^P_t\} = [\ul L^\th_t, \ol L^\th_t]$.
\end{prop}
{ We note that (ii) provides the bounds for the agent's optimal utilities over all possible contracts. As in \reff{Yforward}, later on we will rewrite the agent's continuation utility in a forward form, then such bounds will be crucial for the analysis. }

\begin{rem}
\label{rem-smallesttau}
We remark that the agent's optimal quitting time may not be unique, and among them $\t^{\th; t, \eta}$ is the smallest one. Unfortunately, while the agent's optimal utility is invariant on his choices of the optimal quitting time, the principal's optimal utility does depend on this choice. For technical simplicity in this paper we take the convention that the agent would always take this smallest one $\t^{\th; t, \eta}$, namely he would choose to quit when he is indifferent to quitting or staying.   
\end{rem}

The proof of Proposition \ref{prop-agent1} relies on the following simple lemma, which will be used frequently in the paper. We postpone both proofs to Appendix. 

\begin{lem}
\label{lem-u1domain}
Let Assumption \ref{assum-uAh} hold. Fix $t<T$, $\d>0$,  and recall \reff{Yforward} and \reff{cAy}.

(i) For any $x\le \ol L^\th_t$, there exists a deterministic $\eta^{\th;t,x}\in \cA^P_t$ such that  
\bea
\label{u1domain2}
(\eta^{\th;t,x}, 0) \in \cA^\th_{t,x},\q \int_t^T [C_0 - \eta^{\th;t,x}_s]ds \le \d. 
\eea
More precisely, $\eta^{\th;t,x} := \eta^n$ can be constructed as follows: for some $n$ sufficiently large, 
\bea
\label{etan}
\e_n := {\ol L^\th_t - x\over n} < T-t,\q \eta^n_s := I^f\big(\th; s, f(\th; s, C_0) - n\big) \1_{[t, t+\e_n)} + C_0 \1_{[t+\e_n, T]}.
\eea

(ii) Assume further that Assumption \ref{assum-RCth} holds and $x> \ul L^\th_t$. Then we may and will always construct $\eta^{\th;t,x}$ in a way such that  $(Y, Z, K):= (\cX^{\th; t, x, \eta^{\th;t,x}, 0}, 0, 0)$ solve RBSDE \reff{RBSDE} corresponding to $\eta^{\th; t,x}$. In particular, $\cX^{\th; t, x, \eta^{\th;t,x}, 0}_s > \ol L^\th_s$ for all $s\in [t, T]$ and $\t^{\th; t, \eta^{\th;t,x}}=T$.
\end{lem}

\begin{rem}
\label{rem-u1domain} 
In the rest of the paper we will often use a random version of Lemma \ref{lem-u1domain}, under Assumptions \ref{assum-uAh} and \ref{assum-RCth}. Fix $\th\in \Th$, $(t,x)\in \cA^\th_{t,x}$, $(\eta, Z)\in \cA^\th_{t,x}$, and $\t \in \cT_t$ such that $\t\le \t^{\th; t,\eta}\wedge T_\d$. Then we may construct $\eta^{\th; \t, \cX^{\th; t,x,\eta, Z}_\t}$ in a way such that:
\beaa
(\tilde \eta, \tilde Z) := (\eta, Z)\1_{[t, \t)} + ( \eta^{\th; \t, \cX^{\th; t,x,\eta, Z}_\t}, 0)\1_{[\t, T]} \in  \cA^\th_{t,x}.
\eeaa
Indeed, from the construction \reff{etan} we can easily see that the above $\tilde \eta, \tilde Z$ are $\dbF^t$-progressively measurable, $\tilde Z\in \cA^A_t$, and $\cX^{\th; t,x, \tilde \eta, \tilde Z}_T=0$. Moreover, since $\t\le  \t^{\th; t,\eta}$, then $\cX^{\th; t,x, \tilde \eta, \tilde Z}_\t \in [\ul L^\th_\t, \ol L^\th_\t]$ is bounded. Note further that $T-\t \ge T-T_\d=\d$, then from the construction \reff{etan} we can easily see that $\eta^{\th; \t, \cX^{\th; t,x,\eta, Z}_\t}$ is uniformly bounded, which implies that $\tilde \eta$ has a uniform lower bound, and hence $\tilde \eta\in \cA^P_t$. 
\end{rem}

\subsection{The principal's problem}
Assume the principal hires an agent $\th$ at time $t$ with a contract $\eta \in \cA^P_t$. At time $\t^{\th; t, \eta}$, the agent would quit the job, then the principal will hire a new agent from the market, with possibly a different type $\tilde \th$. In this section, we assume the principal does not want the new agent to quit, by offering only self-enforcing contracts. Then, by the previous section the principal's optimal value for the new contract will be $V^P_0(\tilde \th; \t^{\th; t, \eta})$. The principal will choose $\tilde \th$ to maximize this value, and recall that she will incur a cost $c^P$, so we introduce
\bea
\label{olu0}
\ol u_0(t) := \sup_{\tilde\th\in \Th} V^P_0(\tilde\th; t) - c^P_t.
\eea
{ Clearly $\ol u_0$ inherits the regularity in Proposition \ref{prop-VP0}, and it is bounded. We collect these results in the following simple lemma, and provide a proof in Appendix for completeness. 

\begin{lem}
\label{lem-olu0}
Under  Assumptions \ref{assum-uAh} and \ref{assum-RCth}, $\ol u_0$ is uniformly continuous on $[0, T]$ with a modulus of continuity function $\rho_0 := C\sqrt{\rho}$.  Moreover,  $ \ol u_0$ is bounded  with $\ol u_0(T) = -c^P_T$.
\end{lem} 
}

The principal's problem in this case is:
\bea
\label{VP1}
\left.\ba{c}
\dis V^P_1(\th; t) := \sup_{\eta\in \cA^P_t} \dbE^{\th, Z^{\th;t,\eta}}\Big[ \int_t^{\t^{\th; t, \eta}} \big[ I^h(\th; s,Z^{\th;t,\eta}) - \eta_s\big] ds +\ol u_0(\t^{\th; t, \eta}) \1_{\{\t^{\th; t, \eta}<T\}}\Big],\\
\mbox{subject to}\q V^A_1(\th; t, \eta) \ge R^\th_t.
\ea\right.
\eea
We remark that, here we consider the principal's utility for a given $\th$. In practice, of course, the principal would choose an optimal $\th$ in the beginning, which amounts to computing 
\beaa
V^P_1(t):= \sup_{\th\in \Th} V^P_1(\th; t).
\eeaa

\begin{rem}
\label{rem-adverse}
We emphasize that in our model there are a family of agents whose types are indexed by $\th$. The principal knows the value $\th$ of each agent and can choose the optimal $\th$. This is completely different from the adverse selection problems, which involve only one agent, but the principal does not know the agent's type $\th$ and cannot choose $\th$, cf. \cite{CZ}. 
\end{rem}

Similarly to \reff{Yforward}, we can rewrite RBSDE \reff{RBSDE} in a forward form, which leads to the following dynamic utility function:
\bea
\label{u1}
\left.\ba{c}
\dis u_1(\th; t,x):= \sup_{(\eta, Z)\in \cA^\th_{t,x}} J_1(\th; t, x, \eta, Z),\q (t,x)\in D_\th, \ms\\
\dis  \mbox{where}\q  \t^{\th; t, x, \eta, Z}:= \inf\{s\ge t: \cX^{\th; t, x, \eta, Z}_s = \ul L^\th_s\}\wedge T,\ms\\
\dis J_1(\th; t, x, \eta, Z):= \dbE^{\th,Z}\Big[\int_t^{\t^{\th; t, x, \eta, Z}}\!\!\!\big[ I^h(\th; s,Z_s) - \eta_s\big] ds + \ol u_0(\t^{\th; t, x, \eta, Z}) \1_{\{\t^{\th; t, x,\eta, Z}<T\}}\Big].
\ea\right.
\eea
{ The following result establishes the connection between $V^P_1$ and $u_1$, and thus enables us to use PDE approach for $u_1$ to study $V^P_1$, see the next section.}

\ms
\begin{prop}
\label{prop-V1u1}
Under  Assumptions \ref{assum-uAh} and \ref{assum-RCth}, for any $\th\in \Th$ and $t<T$, we have 
\bea
\label{VP1supu1}
\dis V^P_1(\th; t) = \sup_{R^\th_t \le x\le \ol L^\th_t} u_1(\th; t, x).
\eea
\end{prop}
\proof Let $\tilde J_1(\th; t, \eta)$ denote the expectation  in the right side of \reff{VP1} for a fixed $\eta$. 
For any $\eta\in \cA^P_t$ satisfying the constraint in \reff{VP1}, by Proposition \ref{prop-agent1} we have $x := Y^{\th; t, \eta}_t = V^A_1(\th; t, \eta) \in [R^\th_t, \ol L^\th_t]$. Fix $\d>0$ small and set 
\bea
\label{V1u1taud}
\t_\d := \t^{\th; t, \eta}\wedge T_\d ,~ (\eta^\d, Z^\d)_s := (\eta, Z^{\th; t, \eta})_s \1_{[t, \t_\d )}(s) + (\eta^{\th; \t_\d , Y^{\th; t, \eta}_{\t_\d }}_s, 0) \1_{[\t_\d , T]}(s),
\eea
where $\eta^{\th; \t_\d , Y^{\th; t, \eta}_{\t_\d }}$ is constructed in \reff{etan}. By Remark \ref{rem-u1domain} we have $(\eta^\d , Z^\d )\in \cA^\th_{t, x}$.  By \reff{taueta} we have $Y^{\th; t, \eta}_s> \ul L^\th_s$ and hence $K^{\th; t, \eta}_s =0$ for $s\in [t, \t^{\th; t,\eta})$. Then
\beaa
\cX^{\th; t, x, \eta^\d, Z^\d}_s =  Y^{\th; t, \eta}_s> \ul L^\th_s, ~ s\in [t, \t_\d),\q\mbox{and}\q \t^{\th; t,x, \eta^\d, Z^\d} = \t^{\th; t, \eta} ~\mbox{on}~\{\t^{\th; t, \eta}\le T_\d\}.
\eeaa
Moreover, by Lemma \ref{lem-u1domain} (ii) we have $\t^{\th; t,x, \eta^\d, Z^\d} = T$ on $\{\t^{\th; t, \eta}> T_\d\}$.  Then,
since $\eta^{\th; T_\d , Y^{\th; t, \eta}_{T_\d }}_s$ is $\cF_{T_\d }$-measurable, we have
\beaa
&\dis J_1(\th; t, x, \eta^\d , Z^\d ) = \dbE^{ \th, Z^{\d} }[ \Xi]  = \dbE^{\th,Z^{\th;t,\eta}}[\Xi],\q\mbox{where}\\
&\dis \Xi:= \int_t^{\t_\d} \big[ I^h(\th; s,Z^{\th;t,\eta}_s) - \eta_s\big] ds +  \1_{\{\t^{\th; t, \eta} \le T_\d \}} \ol u_0(\t^{\th; t, \eta})  - \1_{\{\t^{\th; t, \eta} > T_\d \}}  \int_{T_\d }^T \eta^{\th; T_\d , Y^{\th; t, \eta}_{T_\d }}_s ds.
\eeaa
Thus, by omitting the superscripts $^{\th; t, \eta}$, namely denoting $(\t, Y, Z):= (\t^{\th; t, \eta}, Y^{\th; t, \eta}, Z^{\th; t, \eta})$,
 \beaa
&&\dis \tilde J_1(\th; t, \eta) -  u_1(\th; t, x) \le \tilde J_1(\th; t, \eta)  - J_1(\th; t, x, \eta^\d , Z^\d ) \\
&&\dis = \dbE^{\th,Z}\Big[ \1_{\{\t > T_\d \}}  \big[\int_{T_\d }^{\t} [ I^h(\th; s,Z_s) - \eta_s] ds+ \int_{T_\d }^T  \eta^{\th; T_\d , Y_{T_\d }}_s ds + \ol u_0(\t) \1_{\{\t<T\}}\big]\Big]\\
&&\dis \le \dbE^{ \th,Z}\Big[ \int_{T_\d }^T \big[|I^h(\th; s,Z_s)| + |\eta_s|] ds + C \1_{\{T_\d <\t<T\}}\big]\Big] +  C_0\d,
\eeaa
 where the last inequality is due to Lemma \ref{lem-olu0} and the fact $\eta^{\th; T_\d , Y_{T_\d }}_s \le C_0$. Send $\d\to 0$, we obtain $\tilde J_1(\th; t, \eta) \le  u_1(\th; t, x)$. Since $\eta$ is arbitrary, we prove the  "$\le$" part of the proposition.

\ms
To see the opposite direction,  for any $x\in [R^\th_t, \ol L^\th_t]$, $(\eta, Z)\in \cA^\th_{t,x}$, and $\d>0$ small, set
\beaa
\tilde \t_\d:= \t^{\th; t, x, \eta, Z}\wedge T_\d,\q (\tilde \eta^\d, \tilde Z^\d)_s := (\eta, Z)_s \1_{[t, \tilde \t_\d)}(s) + (\eta^{\th; \tilde\t_\d, \cX^{\th; t, x,  \eta, Z}_{\tilde \t_\d}}_s, 0)\1_{[\tilde \t_\d, T]}(s).
\eeaa
Following similar arguments as above, by sending $\d\to 0$ one can show that 
$J_1(\th; t, x, \eta, Z)  \le V^P_1(\th; t)$. 
Then, since $x$ and $(\eta, Z)$ are arbitrary, we obtain the desired  inequality.
\qed

\begin{rem}
\label{rem-outside}
In this section we assume the principal will offer self-enforcing contracts to hire a new agent when the current one quits. There can be other alternatives, for example, the principal may close the business, or run the firm by herself. In those cases, we need to first model the principal's outside option, denoted as $\ol u_0'(t)$, and then study the principal's problem by replacing the $\ol u_0$ in \reff{VP1} with $\ol u_0'$. Provided that $\ol u_0'$ satisfies the technical properties in Lemma \ref{lem-olu0}, the main results in this paper will remain true. 
\end{rem}

\section{The dynamic value function $u_1$}
\label{sect-u1}
\setcounter{equation}{0}
In this section we shall characterize the function $u_1$ as the minimal solution of an HJB equation. { For this purpose, we shall first establish its monotonicity, regularity, and boundary conditions. Their proofs involve very technical constructions. Since they are the main mathematical components of the paper, we report them in the end of this section. However, readers who are not interested in the technical part can skip reading these subsections, which will not affect their understanding on the economic messages of the paper. }

We first show that $u_1$ is decreasing in $x$. { That is, the principal's optimal utility would suffer when the agent's continuation utility $x$ increases. Equivalently, the principal will gain less utility when the agent becomes more expensive, which of course is natural in practice. }

\begin{prop}
\label{prop-u1mon}
Let  Assumptions \ref{assum-uAh} and \ref{assum-RCth} hold. For any  $\th\in \Th$ and $t< T$, the value function $u_1(\theta; t,x)$ is decreasing  in $x\in (\ul L^\th_t, \ol L^\th_t]$. Consequently, 
\bea
\label{VP1u1}
V^P_1(\th; t) =   u_1\big(\th; t, R^\th_t\big).
\eea
\end{prop}

We next show that $u_1$ is uniformly continuous in $(t,x)$. Such regularity is important in its own right, and is crucial for the PDE characterization of $u_1$.

\begin{prop}
\label{prop-u1reg}
Let  Assumptions \ref{assum-uAh} and \ref{assum-RCth} hold. 

(i) For any  $\th\in \Th$ and $t< T$, the value function $u_1(\theta; t,x)$ is uniformly continuous in $x\in (\ul L^\th_t, \ol L^\th_t]$. That is,  for any $\ul L^\th_t< x_1, x_2\le \ol L^\th_t$, we have 
\bea
\label{u1xreg}
\big|u_1(\th; t, x_1) - u_1(\th; t, x_2)\big| \le C\rho_0(|x_1-x_2|),
\eea
where $\rho_0$ is the modulus of continuity function of $\ol u_0$  specified in Lemma \ref{lem-olu0}.

(ii) For any $(t_i, x_i)\in D_\th$, we have 
\bea
\label{u1treg}
|u_1(\th; t_1, x_1) - u_1(\th; t_2, x_2)| \le C\rho_0\big(\rho(|t_1-t_2|) + |x_1-x_2|\big).
\eea
\end{prop}

{ Note that $u_1(\th; t,x)$ involves the hitting time $\t^{\th; t, x, \eta, Z}$, which is in general not continuous in $(t,x)$. This can be circumvented as follows.  Assume we want to estimate $ u_1(\th; t_1, x_1) - u_1(\th; t_2, x_2)$ from above. For any given $(\eta^2, Z^2) \in \cA^\th_{t_2, x_2}$, we want to construct $(\eta^1, Z^1) \in \cA^\th_{t_1, x_1}$ so that $\t^{\th; t_1, x_1, \eta^1, Z^1}$ is close to $\t^{\th; t_2, x_2, \eta^2, Z^2}$ (in addition to other small errors). However, since we have both upper and lower bounds for $\cX$, and since the setting is random, such construction is highly nontrivial, and will rely heavily on Lemma \ref{lem-u1domain}. In particular, we emphasize that the upper bound $\ol L^\th$ and lower bound $\ul L^\th$ play very different roles in the problem, consequently the constructions are quite different for the two cases $x_1<x_2$ and $x_1>x_2$, and even more different for the two cases $t_1<t_2$ and $t_1>t_2$. We refer to Subsection \ref{sect-u1reg} for details. 
}

\begin{rem}
\label{rem-u1mon}
In general the monotonicity property and the continuity of $u_1$ may fail at the boundary point $x=\ul L^\th_t$. For example, for any $x\in (\ul L^\th_t, \ol L^\th_t]$, construct $\eta^n$ as in  \reff{etan},  and send $n\to \infty$, one can easily see that under these contracts the agent never quits and $u_1(\th; t, x) \ge -C_0(T-t)$. However, for $x=\ul L^\th_t$, the agent quits immediately at $t$, and thus $u_1(t, \ul L^\th_t) = \ol u_0(t) = \sup_{\tilde\th} V^P_0(\tilde\th; t)- c^P_t$. Clearly it is possible that $-C_0(T-t) > \ol u_0(t)$ when $c^P_t$ is sufficiently large. 
\end{rem}

Finally, $u_1$ satisfies the following boundary conditions. 
\begin{prop}
\label{prop-u1bound}
Under Assumptions \ref{assum-uAh} and \ref{assum-RCth},  $u_1(\theta; \cd,\cd)$ has boundary conditions:
\bea
\label{u1boundary}
\left.\ba{c}
\dis u_1(\theta; t, \ol L^\th_t)=-C_0(T-t),\q \forall t\in [0, T);\qq u_1(\theta;T-,x)=0,\q \forall  x\in ( \ul L^\th_T,  0);\ms\\
\dis u_1(\th; t, \ul L^\th_t +) \ge \ol u_0(t),\q \forall t\in [0, T).
\ea\right.
\eea
\end{prop}

\bs
We are now ready to characterize $u_1$ through an HJB equation. Introduce
\bea
\label{cLu}
\cL^\th u(t,x) := \pa_t u+\!\!\! \sup_{\eta \leq C_0,z\in \dbR}\!\! \Big[{1\over 2}|z|^2 \partial_{xx} u-\big[f(\th; t, \eta)+\tilde{H}(\th; t, z)\big] \pa_x u +  I^h(\th; t,z)-\eta \Big].
\eea

\begin{thm}
\label{thm-HJB1}
Under Assumptions \ref{assum-uAh} and \ref{assum-RCth}, the value function $u_1(\theta; \cd, \cd)$ is the minimal continuous viscosity solution of the following HJB equation with boundary conditions \reff{u1boundary}:
\bea
\label{HJB1}
\cL^\th u(t,x) =0,\q (t,x)\in D_\th.
\eea
\end{thm}
\proof First, by the regularity in Proposition \ref{prop-u1reg} (i), one may follow rather standard arguments to obtain  the dynamic programming principle: for any $(t, x) \in D_\th$ and $\t \in \cT_t$, 
\bea
\label{u1DPP}
\left.\ba{c}
\dis u_1(\th; t, x) :=  \sup_{(\eta, Z)\in \cA^\th_{t,x}}  \dbE^{\th,Z}\Big[  \int_t^{\t\wedge \t^{\th; t,x,\eta, Z}}[ I^h(\th; r,Z_r) -\eta_r] dr \\
\dis + u_1(\th; \t, \cX^{\th; t, x, \eta, Z}_\t)\1_{\{\t<\t^{\th; t,x,\eta, Z}\}} + \ol u_0(\t^{\th; t,x,\eta, Z}) \1_{\{\t^{\th; t,x,\eta, Z}\le \t, ~ \t^{\th; t,x,\eta, Z}<T\}}  \Big].
\ea\right.
\eea
Then, by \reff{Yforward}, it is clear that $u_1$ is a viscosity solution of \reff{HJB1}. Moreover, by Proposition \ref{prop-u1bound} it satisfies the boundary conditions \reff{u1boundary}. 

We next show that $u_1$ is minimal. Assume $u$ is an arbitrary continuous viscosity solution of \reff{HJB1} satisfying \reff{u1boundary}. In particular, $u(t, \ul L^\th_t) \ge \ol u_0(t)$. For each $n\ge 1$, introduce
\bea
\label{truncate}
\left.\ba{c}
\dis u^{(n)}(t,x) := \sup_{(\eta, Z)\in \cA^P_t \times \cA^A_t, -n \le \eta \le C_0}  \dbE^{ \th,Z}\Big[  \int_t^\t[ I^h(\th; r,Z_r)-\eta_r] dr  + u(\t, \cX^{\th; t,x, \eta, Z}_\t)\Big],\\
\dis\mbox{where}\q \t:=  \ul \t \wedge \ol \t,\q \ul \t:= \t^{\th; t,x, \eta, Z},\q \ol \t:= \inf\{s\ge t: \cX^{\th; t,x, \eta, Z}_s= \ol L^\th_s)\} \wedge T. 
\ea\right.
\eea
Since $\eta$ is uniformly bounded here, one can easily see that $u^{(n)}$ is continuous on $\ol D_\th$, the closure of $D_\th$. Then $u^{(n)}$ is the unique viscosity solution of the truncated HJB equation:
\bea
\label{truncateHJB}
&\dis\!\!\! \pa_t u^{(n)}+\!\!\!\! \sup_{-n\le \eta \leq C_0,z\in \dbR}\!\! \Big[{1\over 2}|z|^2 \partial_{xx} u^{(n)}-\big[f(\th; t, \eta)+\tilde{H}(\th;t,z)\big] \pa_x u^{(n)} +  I^h(\th; t, z)-\eta\Big]=0;\nonumber\\
&\dis \!\!\! u^{(n)}(t, \ol L^\th_t) = -C_0(T-t),\q u^{(n)}(T, x) =0,\q u^{(n)}(t, \ul L^\th_t) = u(t, \ul L^\th_t).
 \eea
Note that $u$ is a viscosity supersolution of  the above HJB equation with the same boundary conditions, then by the comparison principle for \reff{truncateHJB} we have $u\ge u^{(n)}$. 

Now for any $(t, x)\in D_\th$ and  $(\eta, Z)\in \cA^\th_{t,x}$, recall Definition \ref{defn-cA} and let $n \ge C(\eta)$, then $(\eta, Z)$ satisfy the requirements in \reff{truncate}.  Let  $\t, \ul \t, \ol \t$ be defined in \reff{truncate}. Note that
\beaa
u_1(\th; \ol \t, \ol L^\th_{\ol \t}) = -C_0(T-\ol\t)= u(\ol \t, \ol L^\th_{\ol \t}),\q \ol u_0(\ul \t) \le u(\ul\t, \ul L^\th_{\ul\t}) ~\mbox{on}~\{\ul \t<T\}.
\eeaa
Similarly to the DPP \reff{u1DPP}, we have
\beaa
&&\dis J_1(\th; t, x, \eta, Z) \\
&&\dis \le \dbE^{ \th,Z}\Big[  \int_t^{\ol\t\wedge \ul \t} |[ I^h(\th; r,Z_r)-\eta_r] dr +\ol u_0(\ul \t) \1_{\{\ul \t\le \ol \t\}} + u_1(\th; \ol\t, \cX^{\th; t, x, \eta, Z}_{\ol \t})\1_{\{\ul \t>\ol \t\}} \Big]\\
&&\dis =\dbE^{ \th,Z}\Big[  \int_t^\t [ I^h(\th; r,Z_r) -\eta_r] dr + \ol u_0(\ul \t) \1_{\{\ul \t\le \ol \t\}} +u_1(\th; \ol\t, \ol L^\th_{\ol \t})\1_{\{\ul \t>\ol \t\}} \Big]\\
&&\dis \le\dbE^{\th,Z}\Big[  \int_t^\t [I^h(\th; r,Z_r)-\eta_r] dr +u(\ul \t, \ul L^\th_{\ul \t}) \1_{\{\ul \t\le \ol \t\}} + u(\ol \t, \ol L^\th_{\ol \t})\1_{\{\ul \t>\ol \t\}} \Big]\\
&&\dis =\dbE^{\th,Z}\Big[  \int_t^\t [I^h(\th; r,Z_r) -\eta_r] dr +u(\t, \cX^{\th; t, x, \eta, Z}_\t)  \Big] \le u^{(n)}(t,x) \le u(t,x).
\eeaa
Since $(\eta, Z)\in \cA^\th_{t,x}$ is arbitrary, we obtain $u_1(\th; t,x) \le u(t,x)$.
\qed

\begin{rem}
\label{rem-u1jump}
(i) As explained in Remark \ref{rem-u1mon},  in general it is possible that $u_1(\th; t, \ul L^\th_t +) > \ol u_0(t)$. One can show that $u_1(\th;\cd, \cd)$ is the unique viscosity solution to  \reff{HJB1} if we replace the boundary condition at $x=\ul L^\th_t$ with $u_1(\th; t, \ul L^\th_t+)$. However, we are not able to compute $u_1(\th; t, \ul L^\th_t+)$ explicitly, so we characterize $u_1$ as the minimal solution in Theorem \ref{thm-HJB1}.

(ii) The differential operator $\cL^\th$ in \reff{cLu} is the same as the case without constraint on $(\eta, Z)$. The constraint $(\eta, Z)\in \cA^\th_{t,x}$ in \reff{u1} is reflected in the boundary conditions \reff{u1boundary}.
\end{rem}

\subsection{The monotonicity: Proof of Proposition \ref{prop-u1mon}}

First, \reff{VP1u1} follows directly from Proposition \ref{prop-V1u1} and the monotonicity of $u_1$. Without loss of generality we prove the monotonicity only at $t=0$. Fix $\ul L^\th_0< x< x_\d:= x+\d\le \ol L^\th_0$ and assume $\d>0$ is sufficiently small. 
For any $(\eta,Z)\in A^{\theta}_{0,x_\d}$ and $n\ge 1$ such that ${\d\over n}< T$, { our goal  is to construct $(\eta^n, Z^n)\in \cA^\th_{0,x}$ such that 
\bea
\label{etan1}
\!\!\! J_1(\th; 0, x_\d, \eta, Z) - u_1(\th; 0, x) \le J_1(\th; 0, x_\d, \eta, Z) -  J_1(\th; 0, x, \eta^n, Z^n) \to 0,~\mbox{as}~n\to \infty.
\eea
Then by the arbitrariness of $(\eta, Z)$ we obtain $u_1(\th; 0, x_\d) \le u_1(\th; 0, x)$.}

To see this, similarly to \reff{etan} and \reff{V1u1taud} we construct:
\bea
\label{etan2}
\left.\ba{c}
\dis t_n:= {\d\over n},\qq \eta^n_s :=  I^f\big(\th; s, f(\th; s,\eta_s)-n\big), \q Z^n_s := Z_s,\q 0\le s\le t_n;\ms\\
\dis \t_n := \inf\{s\ge 0: \cX^{\th; 0, x, \eta^n, Z^n}_s = \ul L^\th_s\}\wedge t_n,\q \t^{\d} := \t^{\th; 0, x_\d, \eta, Z};\ms\\
\dis (\eta^n_s, Z^n_s) := (\eta^{\th; \t_n, \ul L^\th_{\t_n}}_s, 0) \1_{\{\t_n < t_n\}} +(\eta_s, Z_s) \1_{\{\t_n = t_n\}} ,\q \t_n < s\le T,
\ea\right.
\eea
where $\eta^{\th; \t_n, \ul L^\th_{\t_n}}$ is as in Lemma \ref{lem-u1domain} (with an arbitrary $\d$ in \reff{u1domain2}) and Remark \ref{rem-u1domain} . Denote $\cX^n:= \cX^{\th; 0, x, \eta^n, Z^n}$. It is clear that $(\eta^n, Z^n)\in \cA^P_0\times \cA^A_0$, and by  Lemma \ref{lem-u1domain} we see that $\cX^n_T =0$ on $\{\t_n < t_n\}$.  Note further that $f(\th; s,\eta^n_s)=f(\th; s,\eta_s)-n$, then 
\bea
\label{cXn-d}
\cX^n_s - \cX^{\th; 0, x_\d, \eta, Z}_s  = ns -\d < 0~\mbox{for}~ s\in [0, t_n);\q \mbox{and}\q \cX^n_{t_n} - \cX^{\th; 0, x_\d, \eta, Z}_{t_n}=0.
\eea
Then, since $(\eta^n, Z^n) = (\eta, Z)$ on $[t_n, T]$, we have $\cX^n_T = \cX^{\th; 0, x_\d, \eta, Z}_T =0$ on $\{\t_n=t_n\}$. Therefore, $\cX^n_T=0$ in all the cases, and thus $(\eta^n, Z^n)\in \cA^\th_{0,x}$.

Moreover, note that \reff{cXn-d} implies that 
\beaa
\t^{\th; 0, x, \eta^n, Z^n} = \t_n \1_{\{\t_n < t_n\}} +\t^{\d} \1_{\{\t_n = t_n\}} \le \t^\d, \q  Z^n_s = Z_s,~ 0\le s\le \t^{\th; 0, x, \eta^n, Z^n}.
\eeaa
Then 
\beaa
\dis   J_1(\th; 0, x, \eta^n, Z^n)&=& \dbE^{ \th,Z}\Big[ \1_{\{\t_n < t_n\}}\big[\int_t^{\t_n} [I^h(\th; s,Z_s)- \eta^n_s] ds+ \ol u_0(\t_n) \big]\\
&&\dis  +\1_{\{\t_n = t_n\}}\big[\int_0^{\t^{\d}} [I^h(\th; s,Z_s) - \eta^n_s] ds+ \ol u_0(\t^{\d}) \1_{\{\t^{\d}<T\}}\big]\Big],
\eeaa
and therefore, 
\beaa
&&\dis   J_1(\th; 0, x_\d, \eta, Z) -  J_1(\th; 0, x, \eta^n, Z^n)\\
&&\dis=\dbE^{ \th,Z}\Big[ \int_0^{\t_n} [\eta^n_s-\eta_s]ds +\1_{\{\t_n < t_n\}}\big[\int_{\t_n}^{\t^{\d}} [I^h(\th; s,Z_s) - \eta_s] ds + \ol u_0(\t^{\d}) \1_{\{\t^{\d}<T\}} - \ol u_0(\t_n) \big]\Big]\\
&&\dis \le \dbE^{ \th,Z}\Big[ \1_{\{\t_n < t_n\}}\big[\int_{\t_n}^T [|Z_s| + C_0 - \eta_s] ds +C \big]\Big],
\eeaa
where the last inequality thanks to \reff{HIhproperty}, Lemma \ref{lem-olu0}, and the fact $\eta^n \le \eta$. Moreover, noting that $D_\th$ is bounded, by \reff{Zetaest1}  we have
\beaa
\dbE^{ \th,Z}_{\cF_{\t_n}}\Big[ \int_{\t_n}^T[|Z_s| + C_0 - \eta_s] ds \Big] \le C\sqrt{\ol L^\th_{\t_n} - \cX^\d_{\t_n}} + C(\ol L^\th_{\t_n} - \cX^\d_{\t_n}) \le C.
\eeaa
Then
  \beaa
J_1(\th; 0, x_\d, \eta, Z) - u_1(\th; 0, x) \le  C\dbP^{\th,Z}(\t_n < t_n).
\eeaa
Recall \reff{Yforward} and note that $f<0, \tilde H\le 0$.  Then, by Assumption \ref{assum-RCth} (iii) we have
\beaa
\cX^{\th; 0, x, \eta^n, Z^n}_s - \ul L^\th_s \ge x + \int_0^s Z_r dB_r^{\th,Z} - \ul L^\th_0 - \L_1 t_n,
\q 0\le s\le t_n.
\eeaa
Thus, since $x- \ul L^\th_0 > 0$, for $n$ large enough such that $\L_1 t_n < {x-\ul L^\th_0\over 2}$,
\beaa
\dbP^{ \th,Z}(\t_n < t_n) &\le& \dbP^{ \th,Z}\Big( \inf_{0\le s\le t_n} [\cX^{\th; 0, x, \eta^n, Z^n}_s - \ul L^\th_s] \le 0\Big)\\
&\le&\dbP^{\th,Z}\Big( \inf_{0\le s\le t_n} \int_0^s Z_r dB_r^{\th,Z} \le - {x- \ul L^\th_0\over 2}\Big)\to 0,\q\mbox{as}~n\to \infty.
\eeaa
This implies \reff{etan1}, and hence $ u_1(\th; 0, x_\d)\le u_1(\th; 0, x)$. 
\qed

\subsection{The regularity: Proof of Proposition \ref{prop-u1reg}}
\label{sect-u1reg}
We first need the following technical lemma, whose proof is postponed to Appendix.
\begin{lem}
 \label{lem-Zetaest}
 Let  Assumptions \ref{assum-uAh}  and \ref{assum-RCth} hold. Let $\d>0$ be small, $t\le T_{2\d}$, and $\ul L^\th_t < x \le \ul L^\th_t+\d$. Set $\ul\t_\d := \t^{\th; t,x, C_0, 1} \wedge T_\d$ and 
  \bea
  \label{uletaZd}
\left.\ba{c}
\dis (\ul \eta^{\th; t, x}_s, \ul Z^{\th; t, x}_s) := \1_{[t, {\ul\t}_\d)}(s) (C_0, 1) + \1_{[{\ul \t}_\d, T)}(s)\big(\eta_s^{\th;  \ul\t_\d, \cX^{\th; t,x , C_0, 1}_{{\ul\t}_\d}}, 0\big)\Big].
 \ea\right.
  \eea
Then $(\ul \eta^{\th; t, x}, \ul Z^{\th; t, x})\in \cA^\th_{t,x}$, and 
 \bea
 \label{Zetaest3}
 \ol u_0(t) - J_1(\th; t, x, \ul \eta^{\th; t, x}, \ul Z^{\th; t, x}) \le C\rho_0(\d).
 \eea
 \end{lem}

\bs

\no We now turn to the proof of the proposition.

\no{\bf (i) Proof of \reff{u1xreg}.}  Again without loss of generality we prove it only at $t=0$. Let $x_1 = x$, $x_2= x+\d$ for some $\d>0$ sufficiently small. Then by the monotonicity in Proposition \ref{prop-u1mon} we have $u_1(\th; t, x+\d) - u_1(\th; t, x)\le 0$.  So it suffices to prove the opposite direction. Fix $(\eta, Z)\in \cA^\th_{0,x}$ and, in light of \reff{etan1},  construct
\bea
\label{etaZd}
\left.\ba{c}
\dis (\eta^\d_s, Z^\d_s) := (\eta_s, Z_s)\1_{[0, \t)}(s) +  \1_{[\t, T)}(s) \times\ms\\
\dis \Big[\1_{\{\t = \ol \t\}}(C_0, 0)  + \1_{\{\t =T_{2\d} <\ol\t\wedge \ul \t\}}   \big(\eta^{\th; T_{2\d}, \cX^{\th; 0, x, \eta, Z}_{T_{2\d}}+\d}_t, 0\big)  + \1_{\{\t = \ul \t\}} (\ul \eta^{\th; \t, \ul L^\th_\t + \d}_s, \ul Z^{\th; \t, \ul L^\th_\t + \d}_s\Big];\ms\\
\dis \mbox{where}\q \ol\t:=  \inf\{s\ge 0: \cX^{\th; 0, x, \eta, Z}_s = \ol L^\th_s -\d\};\q \ul\t := \t^{\th; 0, x, \eta, Z};\q \t:= \ol\t \wedge \ul\t \wedge T_{2\d}.
\ea\right.
\eea
 One can verify straightforwardly that $(\eta^\d, Z^\d)\in \cA^\th_{0, x_\d}$, where the lower boundedness of $\eta^\d$ is due to \reff{etan}. Denote
 \beaa
 \cX^x:= \cX^{\th; 0, x, \eta, Z};\q  \cX^\d := \cX^{\th; 0, x_\d, \eta^\d, Z^\d},\q \t^\d := \t^{\th; 0, x_\d, \eta^\d, Z^\d}.
 \eeaa
 Then one can easily check that
\beaa
&\dis \cX^\d_t = \cX_t + \d,~ 0\le t\le \t;\q \cX^\d_\t = \ol L^\th_\t~\mbox{on}~ \{\t=\ol \t\};\q  \cX^\d_\t = \ul L^\th_\t + \d~\mbox{on}~ \{\t=\ul \t\};\\
&\dis \t^\d  > \t ~\mbox{on}~ \{\t = \ul\t\},\q\mbox{and}\q \t^\d = T~\mbox{on}~ \{\t \neq \ul \t\}.
\eeaa
Thus
\bea
\label{I123}
&&J_1(\th; 0, x, \eta, Z) - u_1(\th; 0, x_\d) \le J_1(\th; 0, x, \eta, Z) -J_1(\th; 0, x_\d, \eta^\d, Z^\d) \nonumber\\
&&= \dbE^{\th,Z}\Big[\int_\t^{\ul\t}[I^h(\th; s,Z_s) - \eta_s] ds+ \ol u_0(\ul\t) \1_{\{\ul\t<T\}}\Big]\nonumber\\
&&\q  - \dbE^{\th,Z^\d}\Big[\int_\t^{\t^\d}[I^h(\th; s,Z^\d_s) - \eta^\d_s] ds+ \ol u_0(\t^\d) \1_{\{\t^\d<T\}}\Big]\nonumber\\
&&= \dbE^{\th,Z}\Big[\1_{\{\t=\ol \t\}}  I_1 + \1_{\{\t=T_{2\d} < \ol\t\wedge \ul \t\}} I_2 +   \1_{\{\t=\ul \t\}}  I_3\Big], 
\eea
where
\beaa
\left.\ba{lll}
\dis I_1:= \dbE_{\ol\t}^{\th,Z}\Big[ \int_{\ol\t}^{\ul \t}[I^h(\th; s,Z_s) - \eta_s] ds+\ol u_0(\ul\t) \1_{\{\ul\t<T\}} +C_0(T-\ol\t) \Big];\ms\\
\dis I_2 :=   \dbE_{T_{2\d}}^{ \th,Z}\Big[ \int_{T_{2\d}}^{\ul \t} [I^h(\th; s,Z_s) - \eta_s] ds+ \ol u_0(\ul\t) \1_{\{\ul\t<T\}}  + \int_{T_{2\d}}^T\eta^{\th; T_{2\d}, \cX^\d_{T_{2\d}}}_s ds\big]\Big];\ms\\
\dis I_3 :=\ol u_0(\ul\t) - J_1\big(\th; \ul\t, \ul L^\th_\t+\d, \ul \eta^{\th; \ul\t, \ul L^\th_\t+\d}, \ul Z^{\th; \ul\t, \ul L^\th_\t+\d}\big).
\ea\right.
\eeaa

We now estimate $I_1, I_2, I_3$ separately. First, on $ \{\t=\ol\t\}$, since $\ol u_0$ is bounded, we have
\beaa
I_1&=& \dbE_{\ol\t}^{\th,Z}\Big[ \int_{\ol\t}^{\ul \t}\big[I^h(\th; s,Z_s) +C_0- \eta_s\big] ds\Big]+C \dbP^{\th, Z}_{\ol\t}(\ul\t<T)\\
&\le&\dbE_{\ol\t}^{\th,Z}\Big[ \int_{\ol\t}^T[\L_0|Z_s| +C_0- \eta_s] ds\Big]+C \dbP^{\th, Z}_{\ol\t}(\ul\t<T)\\
&\le& C[\sqrt{\d}+\d] + C \dbP^{\th, Z}_{\ol\t}(\ul\t<T)\le C\sqrt{\d} + C \dbP^{\th, Z}_{\ol\t}(\ul\t<T),
\eeaa
where the second inequality is due to \reff{Zetaest1} and the fact $\ol L^\th_{\ol\t} -\cX^x_{\ol\t} =\d$. 
Moreover, note that $\cX^x_{\ul \t} = \ul L^\th_{\ul \t}\1_{\{\ul\t<T\}}$ and $\cX^x_T=\ol L^\th_T=0$, then by Assumption \ref{assum-RCth} (i) we have
\beaa
0 &\le& \dbE_{\ol\t}^{\th,Z}\Big[ \int_{\ol\t}^{\ul\t} \big[f(\th; r,C_0)-f(\th; r,\eta_r)-\tilde{H}(\th; r,Z_r)\big] dr \Big] \\
&=& \dbE_{\ol\t}^{ \th,Z}\Big[ [\cX^x_{\ul\t} - \cX^x_{\ol\t}] - [\ol L^\th_{\ul \t} - \ol L^\th_{\ol\t}] \Big] \\
&=&  \dbE_{\ol\t}^{ \th,Z}\Big[ \1_{\{\ul\t <T\}} [\ul L^\th_{\ul \t} - \ol L^\th_{\ul \t}] +\d \Big] \le -{1\over \L_0}\dbP^{\th, Z}_{\ol\t}\big(\ul \t<T\big) + \d.
\eeaa
Thus $\dbP^{\th, Z}_{\ol\t}\big(\ul \t<T\big) \le \L_0\d$, and therefore,
\bea
\label{I1est}
I_1 \le C\sqrt{\d} + \L_0\d \le  C\sqrt{\d},\q\mbox{on}\q \{\t=\ol\t\}.
\eea

Next, on $\{\t=T_{2\d} < \ol\t\wedge \ul \t\}$, by \reff{Zetaest1} and noting that $\eta^{\th; T_{2\d}, \cX^\d_{T_{2\d}}}_s\le C_0$, we have
\beaa
I_2 &\le &  \dbE_{T_{2\d}}^{\th,Z}\Big[C \int_{T_{2\d}}^T \big[|Z_s| + |\eta_s| + C_0\big] ds + \ol u_0(\ul\t) \1_{\{\ul\t<T\}}\Big] \le C\sqrt{\d} + C\dbE_{T_{2\d}}^{\th,Z}\big[\ol u_0(\ul\t) \1_{\{\ul\t<T\}}\big].
 \eeaa
Since $\ol u_0(t)$ is uniformly continuous and $\ol u_0(T) = - c^P_T<0$, we have 
 \beaa
\ol  u_0(\ul\t) \le   \ul u_0(\ul\t) - \ol u_0(T) \le \sup_{T_{2\d}\le s\le T} | \ol u_0(s) -\ol u_0(T)| \le \rho_0(2\d)\le 2\rho_0(\d).
 \eeaa
 Here $\rho_0(2\d) \le 2\rho_0(\d)$ is without loss of generality. Thus
 \bea
 \label{I2est}
 I_2 \le C\sqrt{\d} + 2\rho_0(\d) \le C\rho_0(\d),\q\mbox{on}\q \{\t=T_{2\d} < \ol\t\wedge \ul \t\}.
 \eea
 
 Finally, on $ \{\t=\ul\t\}$, by Lemma \ref{lem-Zetaest}  we have $I_3 \le C\rho_0(\d)$. Plug this and \reff{I1est}, \reff{I2est} into \reff{I123}, we get
 \beaa
 J_1(\th; 0, x, \eta, Z) - u_1(\th; 0, x_\d)  \le C\rho_0(\d).
 \eeaa
 Since $(\eta, Z)\in \cA^\th_{0,x}$ is arbitrary, we prove the remaining estimate of \reff{u1xreg} at $t=0$.
 
 \bs
 \no{\bf (ii) Proof of \reff{u1treg}.} For simplicity we assume $t_1 = 0$, $x_1 = x$, $t_2 = \d$, and denote $\D x:= x_2-x_1$.  We assume without loss of generality that $\d>0$ is sufficiently small. 

{\it Step 1.} First, recall Assumption \ref{assum-RCth} (iii) and denote
 \beaa
&\dis \ol t:= \inf\{s\ge 0:  x + \L_1 s  = \ol L^\th_s\}\wedge \d;\q \eta^0_s:= I^f(\th; s, -\L_1)\1_{[0, \ol t]}(s) + C_0 \1_{[(\ol t, \d]}(s),~ 0\le s\le \d;\\
&\dis x_\d := \cX^{\th; 0, x, \eta^0, 0}_\d = x + \L_1 \ol t - \int_{\ol t}^\d f(\th; s, C_0)ds.
 \eeaa
 Then it is clear that
 \beaa
\ul L^\th_s < \cX^{\th; 0, x, \eta^0, 0}_s \le \ol L^\th_s, ~ s\in [0, \d];\q |x_\d - x|\le C\d.
\eeaa
Now for any $(\eta, Z)\in \cA^\th_{\d,x_\d}$,  denote 
\beaa
(\eta^\d_s, Z^\d_s) :=  (\eta^0_s, 0)\1_{[0, \d)}(s) + (\eta_s, Z_s) \1_{[\d, T)}(s).
\eeaa
 One can easily see that  $\cX^{\th; 0, x, \eta^\d, Z^\d}_s = \cX^{\th; \d, x_\d, \eta, Z}_s$ for $s\in [\d, T]$, and $\t^{\th; 0,x,\eta^\d, Z^\d} = \t^{\th; \d, x_\d, \eta, Z}$. Then $(\eta^\d, Z^\d)\in \cA^\th_{0,x}$ and thus
\beaa
J_1(\th; \d, x_\d, \eta, Z) - u_1(\th; 0, x) \le J_1(\th; \d, x_\d, \eta, Z) - J_1(\th; 0,x, \eta^\d, Z^\d) = \int_0^\d \eta^0_sds \le C_0 \d.
\eeaa
By the arbitrariness of $(\eta, Z)\in \cA^\th_{\d,x_\d}$, we obtain $u_1(\th; \d, x_\d) - u_1(\th; 0, x)\le  C_0 \d$. Then, by \reff{u1xreg} and noting that $|x_2 - x_\d|\le |\D x| + C\d$ we have
\beaa
\left.\ba{c}
\dis u_1(\th; \d, x_2) - u_1(\th; 0, x) = u_1(\th; \d, x_2) - u_1(\th; \d, x_\d) + u_1(\th; \d, x_\d) - u_1(\th; 0, x) \\
\dis \le C\rho_0(|\D x|+C \d)  + C_0\d\le C\rho_0\big(|\D x|+\d\big).
\ea\right.
\eeaa
This implies immediately one direction of \reff{u1treg} at $t=0$. 

{\it Step 2.} To see the opposite direction of \reff{u1treg}, we denote
\bea
\label{d'}
\d' := \L_1 \d + 2T\rho(\d).
\eea
We first consider the case that $\ol L^\th_0 - x \le \d'$. Then
\beaa
0\le \ol L^\th_\d - x_2 = \ol L^\th_\d - \ol L^\th_0 + \ol L^\th_0 - x + x- x_2 \le C\rho(\d) + |\D x|.
\eeaa
 Thus, by \reff{u1xreg} we have
\beaa
u_1(\th; 0, x) - u_1(\th; \d, x_2) \le u_1(\th; 0, \ol L^\th_0)  - u_1(\th; \d, \ol L^\th_\d) + C\rho_0\big(\rho(\d)+|\D x|\big).
\eeaa
Moreover, one can easily see  that $u_1(\th; t, \ol L^\th_t) = -C_0(T-t)$. Then
\bea
\label{u1treg2}
u_1(\th; 0, x) - u_1(\th; \d, x_2) \le -C_0\d  + C\rho_0\big(\rho(\d)+|\D x|\big)\le C\rho_0\big(\rho(\d)+|\D x|\big).
\eea

We now consider the case that $x < \ol L^\th_0 - \d'$. Since $x>\ul L^\th_0$, we can easily see that 
\bea
\label{xd'}
x'_\d := x + \d'  \in (\ul L^\th_\d, \ol L^\th_\d).
\eea
 For any $(\eta, Z)\in \cA^\th_{0,x}$, we shall construct in Step 3 below  $(\eta^\d, Z^\d)\in \cA^\th_{\d, x_\d'}$  such that
\bea
\label{J1treg}
J_1(\th; 0, x, \eta, Z) - J_1(\th; \d, x_\d', \eta^\d, Z^\d) \le C\rho_0(\rho(\d)).
\eea
Then $J_1(\th; 0, x, \eta, Z) - u_1(\th; \d, x'_\d) \le C\rho_0(\rho(\d))$. Thus, by the arbitrariness of $(\eta, Z)\in \cA^\th_{0,x}$ we obtain $u_1(\th; 0, x) - u_1(\th; \d, x'_\d) \le C\rho_0(\rho(\d))$. Note further that $ |x'_\d-x_2| \le |\D x|+\d'$, then
\beaa
&\dis u_1(\th; 0, x) - u_1(\th; \d, x_2) = u_1(\th; 0, x) - u_1(\th; \d, x'_\d) + u_1(\th; \d, x'_\d) - u_1(\th; \d, x_2) \\
&\dis \le C\rho_0(\rho(\d)) + C\rho_0(|\D x| + \d') \le  C\rho_0(|\D x| + \rho(\d)).
\eeaa
This, together with \reff{u1treg2}, proves the other direction of \reff{u1treg} in both cases.  

\ms
{\it Step 3.} We now construct $(\eta^\d, Z^\d)\in \cA^\th_{\d, x_\d'}$  satisfying \reff{J1treg}. Denote
\bea
\label{hatcX}
\hat \cX^{\th; \d, x, \eta, Z}_s := x - \int_0^s \big[f(\th; r+\d,\eta_r) + H(\th; r+\d,Z_r)\big] dr + \int_0^s Z_r dB_r,~ 0\le s\le T_\d,
\eea
and introduce the following processes by shifting the space:
\bea
\label{shift}
B'_s := B^{\d}_{s+\d},\q  \eta'_s(B^{\d}_{[\d, s]}) := \eta_{s-\d}(B'_{[0, s-\d]}),\q   Z'_s(B^{\d}_{[\d, s]}) :=  Z_{s-\d}(B'_{[0, s-\d]}).
\eea
Denote  
\bea
\label{shiftcX}
\left.\ba{c}
\dis \cX^x:= \cX^{\th; 0, x, \eta, Z},\q \hat\cX^x:= \hat\cX^{\th; \d, x, \eta, Z},\q  \cX'_t := \hat\cX^x_{t-\d}(B');\ms\\
 \dis \ol\t:=  \inf\big\{t\ge 0: \hat\cX^x_t = \ol L^\th_{t+\d} -\d'\big\};\q \ul\t := \inf\big\{t\ge 0: \cX^x_t = \ul L^\th_t\big\};\q \t:= \ol\t \wedge \ul\t \wedge T_{3\d};\ms\\
\dis  \ol\t':= \ol \t(B') +\d;\q   \ul\t' := \ul \t(B') +\d;\q \t':= \t(B')+\d = \ol\t' \wedge \ul\t' \wedge T_{2\d}.
\ea\right.
\eea
We see that $\Big(B^{\d}_{[\d, T]}, \eta'_{[\d, T]}, Z'_{[\d, T]}, \cX'_{[\d, T]}, \ol\t', \ul\t', \t'\Big)$ and $\Big(B_{[0, T_\d]}, \eta_{[0, T_\d]}, Z_{[0, T_\d]}, \hat \cX^{x}_{[0, T_\d]}, \ol\t+\d, \ul\t+\d, \t+\d\Big)$ have the same joint distribution.  We then construct, similar to \reff{etaZd},
\beaa
&\dis(\eta^\d_s, Z^\d_s) := (\eta'_s,  Z'_s) \1_{[\d, \t')}(s) + \1_{[\t', T)}(s)\times\\
&\dis \Big[ \1_{\{\t' = \ol \t'\}}(C_0, 0) +\1_{\{\t' =T_{2\d} <\ol\t'\wedge \ul \t'\}}   \big(\eta^{\th; T_{2\d}, \cX'_{T_{2\d}}+\d'}_s, 0\big) + \1_{\{\t' = \ul \t'\}} \big(~\ul \eta^{\th; \ul \t', L^\th_{\ul \t'}+\d'}_s, \ul Z^{\th; \ul \t', L^\th_{\ul \t'}+\d'}_s\big) \Big].
\eeaa
Denote $\cX^\d := \cX^{\th; \d, x'_\d, \eta^\d, Z^\d}$ and $\t^\d := \t^{\th; \d, x_\d', \eta^\d, Z^\d}$. It is clear that
\bea
\label{cXd'}
\cX^\d_s = \cX'_s + \d',~ \d\le s\le \t'.
\eea
Moreover, by Assumption \ref{assum-uAh} (iii) and \reff{HIhproperty}, one can easily see that
\bea
\label{hatcXreg}
|\hat\cX^x_s-\cX^x_s|\le 2T\rho(\d),\q 0\le s\le T_\d.
\eea
Then, by Assumption \ref{assum-RCth} (iii),
\beaa
\hat\cX^x_s +\d' \ge \cX^x_s + \d' - 2T \rho(\d) > \ul L^\th_s + \L_1 \d \ge \ul L^\th_{s+\d},\q\mbox{for}\q 0\le s< \ul\t.
\eeaa
By \reff{cXd'} and recalling that $\hat\cX^x$ and $\cX'$ have the same distribution, this implies $\cX^\d_s > \ul L^\th_s$ for $\d\le s < \ul \t'$, a.s. Then one can easily see that $(\eta^\d, Z^\d)\in \cA^\th_{\d, x_\d'}$, and
\beaa
&\dis \cX^\d_{\t'} = \ol L^\th_{\t'} ~\mbox{on}~ \{\t'=\ol \t'\};\q  0\le \cX^\d_{\t'} - \ul L^\th_{\t'-\d} \le 2 \d' ~\mbox{on}~ \{\t'=\ul \t'\};\\
&\dis \t^\d  > \t' ~\mbox{on}~ \{\t' = \ul\t'\},\q\mbox{and}\q \t^\d = T~\mbox{on}~ \{\t' \neq \ul\t'\}.
\eeaa
Note further that
\beaa
 \dbE^{\th,Z^\d}\Big[\int_\d^ {\t'} [I^h(\th; s,Z^\d_s) - \eta^\d_s] ds\Big] &=& \dbE^{\th,Z'}\Big[\int_\d^ {\t(B')+\d} [I^h(\th; s,Z'_s)- \eta'_s] ds\Big]\\
 &=&\dbE^{\th,Z}\Big[\int_0^ {\t}[I^h(\th; s+\d,Z_s) - \eta_s] ds\Big].
 \eeaa
Then
\bea 
\label{I123-t}
J_1(\th; 0, x, \eta, Z) -J_1(\th;  \d, x_\d', \eta^\d, Z^\d)\!\!&=&\!\! \dbE^{\th,Z}\Big[\int_0^{\ul \t} [I^h(\th; s,Z_s) - \eta_s] ds+ \ol u_0(\ul\t) \1_{\{\ul\t<T\}}\Big]\nonumber\\
&&\!\!\!\!\!\!\!\! - \dbE^{\th,Z^\d}\Big[\int_\d^{\t^\d}[I^h(\th; s,Z^\d) - \eta^\d_s] ds+ \ol u_0(\t^\d) \1_{\{\t^\d<T\}}\Big]\nonumber\\
&=& I_0+I_1 + I_2 + I_3,  
\eea
where, by abusing the notations slightly with \reff{I123},  
\beaa
I_0&:=& \dbE^{\th,Z}\Big[\int_0^ {\t}[I^h(\th; s, Z_s)-I^h(\th; s+\d,Z_s)] ds\Big];\\
I_1&:=& \dbE^{\th,Z}\Big[ \1_{\{\t=\ol \t\}} \big[\int_{\ol\t}^{\ul \t}[I^h(\th; s,Z_s) - \eta_s] ds+\ol u_0(\ul\t) \1_{\{\ul\t<T\}}\big]\Big]+ C_0 \dbE^{\th,Z'} \Big[\1_{\{\t'= \ol \t'\}}(T-\ol\t') \Big];\\
I_2 &:=&   \dbE^{\th,Z}\Big[ \1_{\{\t=T_{3\d} < \ol\t\wedge \ul \t\}}\big[ \int_{T_{3\d}}^{\ul \t}[I^h(\th; s,Z_s) - \eta_s] ds+ \ol u_0(\ul\t) \1_{\{\ul\t<T\}}\big]\Big]  \\
&&+ \dbE^{\th,Z'} \Big[\1_{\{\t'=T_{2\d} < \ol\t'\wedge \ul \t'\}} \int_{T_{2\d}}^T  \eta^{\th; T_{2\d}, \cX^\d_{T_{2\d}}}_s ds\big]\Big];\\
I_3 &:=& \dbE^{\th,Z}\Big[\1_{\{\t=\ul \t\}} \ol u_0(\ul\t)\Big] - \dbE^{\th,Z'}\Big[  \1_{\{\t'=\ul \t'\}} J_1\big(\th; \ul\t', \cX^\d_{\ul \t'}, \ul \eta^{\th; \ul\t', \cX^\d_{\ul \t'}}, \ul Z^{\th; \ul\t', \cX^\d_{\ul \t'}}\big)\Big].
\eeaa

We now estimate $I_i$. First, by the regularity of $h$ we have $|I_0| \le C\rho(\d)$. Next, note that
\beaa
&&\dis \dbE^{\th,Z'} \Big[\1_{\{\t'= \ol \t'\}}(T-\ol\t') \Big]= \dbE^{\th,Z} \Big[\1_{\{\t= \ol \t\}}[ T-\t-\d] \Big]\le \dbE^{\th,Z} \Big[\1_{\{\t= \ol \t\}}[T-\t] \Big];\\
&&\dis \ol L^\th_{\ol \t} - \cX^x_{\ol\t} \le \ol L^\th_{\ol \t} - \hat \cX^x_{\ol \t} + \d' = \ol L^\th_{\ol \t} - \ol L^\th_{\ol \t+\d} + 2\d'\le 2\d'.
\eeaa
Then, similarly to \reff{I1est} and  \reff{I2est}, we have 
\beaa
I_1 \le C\sqrt{\d'}\le C\sqrt{\rho(\d)}\q\mbox{and}\q I_2\le C\rho_0(\d') \le C\rho_0(\rho(\d)).
\eeaa
 Finally,  when $\ul\t=\t\le T_{3\d}$, by \reff{shiftcX} and \reff{hatcXreg} we have
\beaa
\hat \cX^x_{\ul\t} - \ul L^\th_{\ul \t} = \hat \cX^x_{\ul\t} - \cX^x_{\ul\t} \le C\rho(\d).
\eeaa
Again due the the related identical distribution, we have $\cX'_{\ul\t'} - \ul L^\th_{\ul \t'-\d} \le C\rho(\d)$ a.s. on $\{\ul\t=\ul\t'\}$. Then, by \reff{cXd'} and the regularity of $\ul L^\th$,
\beaa
\cX^\d_{\ul\t'} - \ul L^\th_{\ul \t'} \le \cX'_{\ul \t'} + \d' -  \ul L^\th_{\ul \t'-\d}  + \L_1 \d\le C\rho(\d),\q \mbox{a.s.}
\eeaa
Now it follows from Lemma \ref{lem-Zetaest} (iii) that
\beaa
I_3 &\le& \dbE^{\th,Z}\Big[\1_{\{\t=\ul \t\}} \ol u_0(\ul\t)\Big] - \dbE^{\th,Z'}\Big[  \1_{\{\t'=\ul \t'\}} \ol u_0(\ul \t')\Big] + C\rho_0(\rho(\d)) \\
&=&\dbE^{\th,Z}\Big[\1_{\{\t=\ul \t\}} \ol u_0(\ul\t)\Big] - \dbE^{\th,Z}\Big[  \1_{\{\t=\ul \t\}} \ol u_0(\ul \t +\d)\Big] + C\rho_0(\rho(\d)) \le C\rho_0(\rho(\d)).
\eeaa
 Plug all these estimates into \reff{I123-t}, we  obtain \reff{J1treg}.
\qed

\subsection{The boundary conditions: Proof of Proposition \ref{prop-u1bound}}
Since $\cA^\th_{t, \ol L^\th_t} = \{(C_0,0)\}$, the equality $u_1(\theta; t, \ol L^\th_t)=-C_0(T-t)$ is obvious, and actually has already been used in the proof of Proposition \ref{prop-u1reg} (ii). Next, by using  Lemma \ref{lem-u1domain},
\beaa
u_1(\th; t, x) \ge J_1(\th; t, x, \eta^{\th; t, x}, 0) = -\int_t^T  \eta^{\th; t, x}_s ds \ge - C_0(T-t).
\eeaa
Then $\liminf_{t\uparrow T} u_1(\theta;t,x)\ge 0$. On the other hand, for any $ \ul L^\th_T < x<  0$, any $t<T$ large enough such that $\ul L^\th_t < x$, and any $(\eta, Z)\in \cA^\th_{t,x}$, by \reff{Zetaest1} and Lemma \ref{lem-olu0} we have
\beaa
&&\dis J_1(\th; t, x, \eta, Z) = \dbE^{\th,Z}\Big[\int_t^T[I^h(\th; s,Z_s)-\eta_s] ds+ \ol u_0(\t^{\th; t, x,\eta, Z})\1_{\{\t^{\th; t, x,\eta, Z}<T\}} \Big]\\
&&\dis \le \dbE^{\th,Z}\Big[\int_t^T\big[C|Z_s|+ |\eta_s|\big] ds+ \1_{\{\t^{\th; t, x,\eta, Z}<T\}} \big[\ol u_0(T) + \sup_{t\le s\le T} |\ol u_0(s)-\ol u_0(T)| \big]\Big]\\
&&\le C\sqrt{T-t} + [\rho_0(T-t)- c^P_T] \dbP\big(\t^{\th; t, x,\eta, Z}<T\big) \le C\sqrt{T-t} + \rho_0(T-t).
\eeaa
By the arbitrariness of $(\eta, Z)$, we have $u_1(\theta;t,x) \le C\sqrt{T-t} + \rho_0(T-t)$. Send $t\uparrow T$, together with $\liminf_{t\uparrow T} u_1(\theta;t,x)\ge 0$, we obtain $u_1(\theta;T-,x)=0$. 

Finally, fix $t<T$. For any $\d>0$ small, by \reff{Zetaest3} we have $\ol u_0(t) \le u_1(\th; t, x) + C\rho_0(\d)$ for all $x\in (\ul L^\th_t, \ul L^\th_t+\d]$.  
Send $\d\downarrow 0$, since $u_1$ is decreasing, we obtain  $u_1(\th; t, \ul L^\th_t +) \ge \ol u_0(t)$.
\qed

\bs

\section{The case with at most $n$ quittings}
\label{sect-nquit}
\setcounter{equation}{0}
In this section we extend the results in the previous section to the case that the principal allows for at most $n$ quittings. That is, if there have been $n$ agents quitting the job, the principal will only offer a self-enforcing contract so that  the $(n+1)$-th agent will not quit anymore. We remark that, for each agent, the agent's problem remains the same as in the previous section, in particular, we do not keep track of the agent after he quits the job. Then, given the setting in the beginning of Section \ref{sect-1quit}, in particular under \reff{ulLtht} and Assumption \ref{assum-RCth}, we introduce the principal's dynamic utility function recursively as in \reff{u1}: 
\bea
\label{un}
\left.\ba{c}
\dis u_n(\th; t,x) :=\!\!\!\! \sup_{(\eta, Z)\in \cA^\th_{t,x}} \!\!\!\!\dbE^{ \th,Z}\Big[\int_t^{\t^{\th; t, x, \eta, Z}} \!\!\!\!\!\!\!\!\!\! [I^h(\th; s,Z_s)- \eta_s] ds+  \ol u_{n-1}(\t^{\th; t, x, \eta, Z}) \1_{\{\t^{\th; t, x,\eta, Z}<T\}}\Big];\ms\\
\dis \ol u_n(t) := \sup_{\th\in \Th} u_n(\th; t, R^\th_t) - c^P_t,
\ea\right.
\eea
for $n\ge 1$ and  $(t, x)\in D_\th$, and as in \reff{VP1} the principal's problem is:
\bea
\label{VPn}
\left.\ba{c}
\dis V^P_n(\th; t) := \sup_{\eta\in \cA^P_t} \dbE^{ \th, Z^{\th; t, \eta}}\Big[ \int_t^{\t^{\th; t, \eta}} \big[ I^h(\th; s, Z^{\th; t, \eta}_s) - \eta_s\big] ds  +\ol u_{n-1}(\t^{\th; t, \eta}) \1_{\{\t^{\th; t, \eta}<T\}}\Big],\\
\mbox{subject to}\q V^A_1(\th; t, \eta) \ge R^\th_t.
\ea\right.
\eea

By the same arguments as in the previous section, the following results are obvious.

\begin{thm}
\label{thm-HJBn}
(i) $\ol u_n$ is bounded and uniformly continuous on $[0, T]$ with certain modulus of continuity function $\rho_n$, and $\ol u_n(T) = -c^P_T$.

(ii) $u_n$ is decreasing in $x$. In particular, $V^P_n(\th; t) = u_n(\th; t, R^\th_t)$.

(iii) For any  $(t_i, x_i)\in D_\th$, $i=1,2$, with $\D t:= t_2-t_1$, $\D x:= x_2-x_1$,  we have
\bea
\label{untreg}
|u_n(\th; t_1, x_1) - u_n(\th; t_2, x_2)| \le C\rho_{n-1}(\rho(|\D t|) + |\D x|).
\eea

(iv) The value function $u_n(\theta; \cd, \cd)$ is the minimal continuous viscosity solution of the HJB equation \reff{HJB1} with the following boundary conditions:
\bea
\label{unboundary}
\left.\ba{c}
\dis u_n(\theta; t, \ol L^\th_t)=-C_0(T-t),\q \forall t\in [0, T);\qq u_n(\theta;T-,x)=0,\q \forall  x\in (\ul L^\th_T,  0);\ms\\
\dis u_n(\th; t, \ul L^\th_t +) \ge \ol u_{n-1}(t),\q \forall t\in [0, T).
\ea\right.
\eea
\end{thm}

Moreover, we have the following monotonicity result.
\begin{prop}
\label{prop-unmon}
Under Assumptions \ref{assum-uAh} and \ref{assum-RCth},  $u_n$ and $V^{P}_n$ are increasing in $n$.
\end{prop}
\proof For any $\th\in \Th$, $t<T$, and $\eta\in \cA^P_t$ satisfying the constraint in \reff{VPt}, by setting $\t=T$ in \reff{VA1} we have
$V^A_1(\th; t, \eta)\ge Y^{\th; \eta}_t \ge R^{\th}_t$, where $Y^{\th; \eta}$ solves the BSDE \reff{BSDE}. Moreover, by \reff{RBSDE} and \reff{taueta} clearly  $(Y^{\th; t, \eta}, Z^{\th; t, \eta}, K^{\th; t, \eta})= (Y^{\th; \eta}, Z^{\th; \eta}, 0)$  and $\t^{\th; t, \eta}=T$. Then by \reff{VPt} and \reff{VP1} we see that $V^P_1(\th; t) \ge J_P(\th; t, \eta)$, and thus $V^P_1(\th; t) \ge V^P_0(\th; t)$. Now by \reff{VP1u1}, \reff{VP1}, and \reff{VPn} we see that $u_2 \ge u_1$. Then by \reff{VPn} again we may prove recursively that $u_n$ is increasing in $n$, which implies further that $V^{P}_n$ is also increasing.
\qed

\begin{rem}
\label{rem-unmon}
{ (i) One might think that $V^P_n$ should be decreasing in $n$ at the first thought, since, for a larger $n$, the principal is giving agents more opportunity to quit and in the meantime the principal would incur a cost $c^P$. This increasing property is due to our setting that the principal offers self-enforcing contracts after $n$ quittings. Indeed, in the formulation of $V^P_n$, the principal can still allow the $(n+1)$-th agent to quit, but since the contract is self-enforcing, he has incentive to quit. While for the problem $V^P_{n+1}$, the principal has more freedom to choose contracts for the $(n+1)$-th agent, and therefore her optimal utility may increase. Such monotonicity may not hold true if we replace the principal's outside option $\ol u_0(t)$ with another $\ol u_0'(t)$ as in Remark \ref{rem-outside}. In fact, the opposite monotonicity is also possible, see Proposition \ref{prop-uSnmon} below. }

\ms
(ii) The above monotonicity holds true even when there is only one type of agent, namely $\Th=\{\th_0\}$ is a singleton, as in the standard literature. Again this proposition implies that, it is indeed not optimal for the principal to restrict to self-enforcing contracts. She would rather allow the agents to quit and then hire a new one when the current one quits the job. 
\end{rem}

\subsection{An alternative representation formula for $u_n$}
In this subsection we provide an alternative representation formula for $u_n$, which is not recursive and will be used in the next section. To illustrate the idea, consider $u_2(\th_0; t, x)$ for some $\th_0\in \Th$ and $(t, x)\in D_{\th_0}$. First, there exists an (approximately) optimal  $(\eta^0, Z^0)\in \cA^{\th_0}_{t, x}$. Denote $\t_1:= \t^{\th_0; t, x, \eta^0, Z^0}$. When $\{\t_1 < T\}$, by \reff{un}  we may choose an (approximately) optimal $\th_1$, which is an $\cF_{\t_1}$-measurable $\Th$-valued random variable,  such that 
\beaa
\ol u_1(\t_1) \approx u_1(\th_1; \t_1, R^{\th_1}_{\t_1}) - c^P_{\t_1}.
\eeaa
We may continue to find (approximately) optimal $(\eta^1, Z^1)\in \cA^{\th_1}_{\t_1, R^{\th_1}_{\t_1}}$ for $u_1(\th_1; \t_1, R^{\th_1}_{\t_1})$. Denote $\t_2:= \t^{\th_1; \t_1, R^{\th_1}_{\t_1}, \eta^1, Z^1}$, then formally we may rewrite \reff{un} as
\bea
\label{u1-alternative}
\left.\ba{c}
\dis u_1(\th_0; t, x) = \sup_{(\eta^0, Z^0)\in \cA^{\th_0}_{t,x}} \sup_{\th_1} \sup_{(\eta^1, Z^1) \in \cA^{\th_1}_{\t_1, R^{\th_1}_{\t_1}}} \dbE^{\th_0, Z^0}\Big[\int_t^{\t_1} [I^h(\th_0; s, Z^0_s) - \eta^0_s] ds \\
\dis +  \dbE_{\cF_{\t_1}}^{\th_1,Z^1} \big[\int_{\t_1}^{\t_2} [I^h(\th_1; s,Z^1_s)- \eta^1_s] ds + \ol u_0(\t_2)\1_{\{\t_2 < T\}} \big] - c^P_{\t_1}\1_{\{\t_1<T\}} \Big].
\ea\right.
\eea 
Moreover, denote
\bea
\label{zeta}
\zeta_s:= \th_0 \1_{[t, \t_1)}(s) + \th_1 \1_{[\t_1, T)}(s),\q (\eta_s, Z_s) :=  (\eta^0_s, Z^0_s) \1_{[t, \t_1)}(s) + (\eta^1_s, Z^1_s) \1_{[\t_1, T)}(s),
\eea
and let $\dbL^0(\dbF^t; \Th)$ denote the set of $\dbF^t$-progressively measurable $\Th$-valued processes. Then the expectation in \reff{u1-alternative} can be rewritten as:
\bea
\label{EzetaZ}
\dbE^{\zeta, Z}\Big[\int_t^{\t_2} [I^h(\zeta_s; s, Z_s) - \eta_s] ds  - c^P_{\t_1}\1_{\{\t_1<T\}} + \ol u_0(\t_2)\1_{\{\t_2 < T\}}\Big].
\eea
Recall \reff{BMO} it is clear that $Z\in \cA^A_t$. However, we remark that the lower bound of $\eta^1$ may not be uniform, and thus in general the above $\eta\notin \cA^P_t$.  

We now present rigorously the  alternative representation formula for $u_n$. Fix $\th\in \Th$ and $(t, x) \in D_{\th}$, let $\cA^{\infty, \th}_{t,x}$ denote the set of $(\eta, Z, \zeta)\in \cA^P_t \times \cA^A_t \times \dbL^0(\dbF^t; \Th)$ such that there exist associated  $(\{\t_i\}_{i\ge 0}, \cX)$ satisfying:
\bea
\label{cAinfty}
\left.\ba{c}
\dis  \t_0 := t,\q  \cX_{\t_0} := x,\q \zeta_{\t_0} = \th,\q   \mbox{and for}~ i\ge 0,\ms\\
\dis  \t_{i+1}:= \inf\{s\ge \t_i:   \cX^{\zeta_{\t_i}; \t_i, \cX_{\t_i}, \eta, Z}_s = \ul L^{\zeta_{\t_i}}_s\} \wedge T;\ms\\
\dis  \zeta_s =  \zeta_{\t_i}, \q \cX_s := \cX^{\zeta_{\t_i}; \t_i, \cX_{\t_i}, \eta, Z}_s, ~ \t_i \le s < \t_{i+1};\ms\\
\dis \cX_{\t_{i+1}} := R^{\zeta_{\t_{i+1}}}_{\t_{i+1}}\1_{\{\t_{i+1}< T\}} +\cX^{\zeta_{\t_i}; \t_i, \cX_{\t_i}, \eta, Z}_T\1_{\{\t_{i+1}=T\}}.
\ea\right.
\eea
Here $\t_i$ stands for the quitting time of the $i$-th agent, and when he quits, the principal would hire a new agent with type $\zeta_{\t_i}$. We note that $\cX$ have jumps at $\t_i$.

\begin{prop}
\label{prop-unalternative}
Let Assumptions \ref{assum-uAh} and \ref{assum-RCth} hold. Then 
\bea
\label{un2}
\left.\ba{c}
\dis u_n(\th; t, x) = \sup_{(\eta, Z, \zeta)\in \cA^{\infty, \th}_{t,x}} J_n(t, x, \eta, Z, \zeta),\q\mbox{subject to}\q \cX_T=0~\mbox{on}~\{\t_n=T\}, \q\mbox{where}\\
\dis  J_n(t, x, \eta, Z, \zeta):= \dbE^{\zeta,Z}\Big[\int_t^{\t_n} [I^h(\zeta_s; s,Z_s) - \eta_s] ds  -\sum_{i=1}^{n-1} c^P_{\t_i}\1_{\{\t_i<T\}} + \ol u_0(\t_n)\1_{\{\t_n<T\}}\Big].
\ea\right.
\eea
\end{prop}
\no The proof is purely technical, and is postponed to Appendix.

\section{The general case with arbitrary number of quittings}
\label{sect-general}
\setcounter{equation}{0}

In light of \reff{un2}, we first introduce the principal's dynamic value function when there is no constraint on the number of agent's quittings: for any $\th\in \Th$ and $(t, x) \in D_{\th}$, 
\bea
\label{uinfty}
\left.\ba{c}
\dis u_\infty(\th; t, x) = \sup_{(\eta, Z, \zeta)\in \cA^{\infty, \th}_{t,x}} J_\infty(t, x, \eta, Z, \zeta), \q\mbox{subject to}\q \lim_{i\to\infty}\t_i=T, \cX_{T-}=0,\\
\dis \mbox{where}\q J_\infty(t, x, \eta, Z, \zeta):= \dbE^{\zeta,Z}\Big[\int_t^T [I^h(\zeta_s; s,Z_s) - \eta_s] ds  -\sum_{i=1}^\infty c^P_{\t_i}\1_{\{\t_i<T\}} \Big].
\ea\right.
\eea
Then, if the principal hires an agent $\th$ at time $t$ and allows the agent as well as all future agents to quit as they want, the principal's value at time $t$ will be:
\bea
\label{VPinfty}
V^P_\infty(\th, t) := \sup_{R^{\th}_t \le x\le \ol L^{\th}_t} u_\infty(\th; t, x),\q V^P_\infty(t) := \sup_{\th\in \Th} V^P_\infty(\th, t).
\eea
Here $V^P_\infty(t)$ will be the principal's real optimal utility if he hires an agent at time $t$, by hiring the best agent for her. 

We remark that,  in general it may look possible that agents quit too frequently and there are already infinitely many agents quitting at some time before $T$, namely $\t_\infty:= \lim_{i\to \infty} \t_i < T$.  { However, due to our assumption $c^\th_t \ge {1\over \L_0}>0$, this can never be the case. Roughly speaking, each new agent is hired with continuation utility $R^\th$, and he will remain in the position until his continuation utility hits the lower barrier $\ul L^\th = R^\th- c^\th \le R^\th - {1\over \L_0}$. Then, by the desired regularity of $R^\th$ and $\ul L^\th$, it will take a while for the agent to reach the lower barrier. That is, the time period $\t_i-\t_{i-1}$ during which the $i$-th agent works for the principal will have a lower bound (in certain random sense), consequently, within a finite horizon $[0, T]$ there can be only finitely many agents.}  Indeed, the following result shows that we shall have $\t_n = T$, for $n$ large enough, a.s. for all $(\eta, Z, \zeta)\in \cA^{\infty, \th}_{t,x}$, and thus the constraint  $\dis\lim_{i\to\infty}\t_i=T$ in \reff{uinfty} is actually redundant. In particular, there is no self-enforcing contract involved anymore in \reff{uinfty} or \reff{VPinfty}. The proof is again postponed to Appendix.

\begin{prop}
\label{prop-taufinite}
Let Assumptions \ref{assum-uAh} and \ref{assum-RCth} hold. Fix arbitrary $\th\in \Th$ and $(t,x)\in D_{\th}$.

(i) For any  $(\eta, Z, \zeta)\in \cA^{\infty,\th}_{t,x}$ with corresponding $(\{\t_i\}_{i\ge 0}, \cX)$ defined in \reff{cAinfty}, we have 
\bea
\label{taun}
\dbP^{\zeta,Z}(\t_n < T) \le {C\over n},~ \forall n\ge 1, \q\mbox{and consequently}, \q \dbP^{\zeta,Z}\big(\t_n < T, ~\mbox{for all $n$}\big) =0.
\eea

{ (ii)  $u_n$ converges to $u_\infty$ uniformly:}
\bea
\label{unconv}
|u_n(\th; t, x) - u_\infty(\th; t, x)|\le {C\over n}.
\eea
\end{prop}

The following result is a direct consequence of Theorem \ref{thm-HJBn} and \reff{unconv}. Define
\bea
\label{oluinfty}
\ol u_\infty(t) := \sup_{\th\in \Th} u_\infty(\th; t, R^\th_t) - c^P_t.
\eea

\begin{prop}
\label{prop-conv}
Let Assumptions \ref{assum-uAh} and \ref{assum-RCth} hold. Then,
 
(i) $u_\infty$ is decreasing in $x$. In particular, 
\beaa
V^P_\infty(\th; t) = u_\infty(\th; t, R^\th_t),\q\mbox{and thus}\q \ol u_\infty(t)  = V^P_\infty(t) - c^P_t.
\eeaa

(ii) $u_\infty$ is uniformly continuous in $D_\th$, uniformly in $\th$. \

(iii) $\ol u_\infty$ is bounded and uniformly continuous on $[0, T]$ with $\ol u_\infty(T) = -c^P_T$.

(iv)  $u_\infty(\theta; t,x)$ has boundary conditions:
\bea
\label{uinftyboundary}
\left.\ba{c}
\dis u_\infty(\theta; t, \ol L^\th_t)=-C_0(T-t),\q \forall t\in [0, T);\qq u_\infty(\theta;T-,x)=0,\q  \forall  x\in (\ul L^\th_T,  0);\ms\\
\dis u_\infty(\th; t, \ul L^\th_t +) \ge \ol u_\infty(t),\q \forall t\in [0, T).
\ea\right.
\eea
\end{prop}

\ms
We now extend Theorem \ref{thm-HJBn} (iv).
\begin{thm}
\label{thm-HJBinfty}
The value function $\{u_\infty(\theta; \cd, \cd)\}_{\th\in \Th}$ is the minimal continuous viscosity solution of the HJB system \reff{HJB1} for all $\th\in \Th$ with boundary conditions \reff{uinftyboundary}.
\end{thm}

\begin{rem}
\label{rem-HJBinfty}
In Theorem \ref{thm-HJBn},  the boundary condition $\ol u_{n-1}$ is given (recursively), and then $u_n(\th; \cd,\cd)$ can be solved through \reff{HJB1} for fixed $\th$. The situation here is quite different.  By \reff{oluinfty}, the boundary condition $\ol u_\infty$ relies on the solution $u_\infty$ itself, and it involves all $\th$. So here the equations for different $\th$ are interacting through the boundary condition and thus we are considering an infinite dimensional system of HJB equations.   
\end{rem}

\no{\bf Proof of Theorem \ref{thm-HJBinfty}.} By extending \reff{u1DPP} to $u_n$ and sending $n\to \infty$, it follows from the uniform convergence of $u_n$ that $u_\infty$ satisfies the dynamic programming principle:
\bea
\label{uinftyDPP}
\left.\ba{c}
\dis u_\infty(\th; t, x) :=  \sup_{(\eta, Z)\in \cA^\th_{t,x}}  \dbE^{\th,Z}\Big[  \int_t^{\t\wedge \t^{\th; t,x,\eta, Z}} [ I^h(\th; r,Z_r) -\eta_r] dr \\
\dis + u_\infty(\th; \t, \cX^{\th; t, x, \eta, Z}_\t)\1_{\{\t<\t^{\th; t,x,\eta, Z}\}} + \ol u_\infty(\t^{\th; t,x,\eta, Z}) \1_{\{\t^{\th; t,x,\eta, Z}\le \t,  \t^{\th; t,x,\eta, Z}<T\}}\Big],
\ea\right.
\eea
 for any $\th\in \Th$, $(t, x) \in D_\th$, and $\t \in \cT_t$.
Then it is straightforward to show that $\{u_\infty(\theta; \cd, \cd)\}_{\th\in \Th}$ is a continuous viscosity solution of the HJB system \reff{HJB1} for all $\th\in \Th$ with boundary conditions \reff{uinftyboundary}. Thus it remains to prove the minimal property. 

Let $\{u(\theta; \cd, \cd)\}_{\th\in \Th}$ be an arbitrary continuous viscosity solution of the HJB system \reff{HJB1} for all $\th\in \Th$ with boundary conditions \reff{uinftyboundary}. We emphasize that here the boundary condition at $\ul L^\th_t$ is not $\ol u_\infty$, but rather 
\beaa
 \ol u(t) := \sup_{\th\in \Th} u(\th; t, R^\th_t) - c^P_t.
\eeaa
Fix $(\th, t, x)$ and $(\eta, Z, \zeta)\in \cA^{\infty,\th}_{t,x}$ with corresponding $(\{\t_i\}_{i\ge 0}, \cX)$ satisfying $\cX_{T-}=0$. Given $\ol u$, recall \reff{u1} and introduce
\beaa
\left.\ba{c}
\dis \tilde u_1(\th; t,x) :=\!\!\! \sup_{(\eta^0, Z^0)\in \cA^\th_{t,x}} \!\!\! \dbE^{ \th,Z^0}\Big[\int_t^{\t^{\th; t, x, \eta^0, Z^0}} \!\!\!\!\!\!\! [I^h(\th; s,Z^0_s) - \eta^0_s] ds+ \ol u(\t^{\th; t, x, \eta^0, Z^0}) \1_{\{\t^{\th; t, x,\eta^0, Z^0}<T\}}\Big].
\ea\right.
\eeaa
Note that $u(\th; \cd,\cd)$ is a viscosity solution to \reff{HJB1} with boundary condition $\ol u$. Then by Theorem \ref{thm-HJB1} we see that $\tilde u_1(\th; t, x) \le u(\th; t, x)$. Thus for any $(\eta^0, Z^0)\in  \cA^\th_{t,x}$ we have
\beaa
\dbE^{\th,Z^0}\Big[\int_t^{\t^{\th; t, x, \eta^0, Z^0}}[I^h(\th; s,Z^0_s) - \eta^0_s] ds+  \ol u(\t^{\th; t, x, \eta^0, Z^0}) \1_{\{\t^{\th; t, x,\eta^0, Z^0}<T\}}\Big]  \le u(\th; t, x).
\eeaa
Note that in general $(\eta, Z)\notin\cA^\th_{t,x}$, and thus at above one cannot set $(\eta^0, Z^0)=(\eta, Z)$. However, following the arguments in Proposition \ref{prop-unalternative} one can still derive from above  that
\beaa
\dbE^{\th,Z}\Big[\int_t^{\t_1}[I^h(\th; s,Z_s) - \eta_s] ds+ \ol u(\t_1) \1_{\{\t_1<T\}}\Big] \le u(\th; t, x). 
\eeaa
Moreover, by the definition of $\ol u$, we have $u(\zeta_{\t_1}; \t_1, R^{\zeta_{\t_1}}_{\t_1}) - c^P_{\t_1} \le \ol u(\t_1)$ on $\{\t_1 < T\}$, then
\beaa
 \dbE^{ \th,Z}\Big[\int_t^{\t_1}[ I^h(\th; s,Z_s) - \eta_s] ds+ \big[u(\zeta_{\t_1}; \t_1, R^{\zeta_{\t_1}}_{\t_1}) - c^P_{\t_1}\big] \1_{\{\t_1<T\}}\Big] \le u(\th; t, x).
\eeaa
Similarly, on $\{\t_1 < T\}$ we have
\beaa
 \dbE_{\cF_{\t_1}}^{\zeta_{\t_1},Z}\Big[\int_{\t_1}^{\t_2}[I^h(\zeta_{\t_1}; s,Z_s) - \eta_s] ds+ \ol u(\t_2) \1_{\{\t_2<T\}}\Big] \le u(\zeta_{\t_1}; \t_1, R^{\zeta_{\t_1}}_{\t_1}).
\eeaa
Then
\beaa
\dbE^{\zeta; \cd,Z}\Big[\int_t^{\t_2}[I^h(\zeta_s; s,Z_s) - \eta_s] ds+ \ol u(\t_2)\1_{\{\t_2<T\}}  - c^P_{\t_1} \1_{\{\t_1<T\}}\Big] \le u(\th; t, x). 
\eeaa
Repeating the arguments we have, for any $n\ge 1$,
\beaa
\dbE^{\zeta,Z}\Big[\int_t^{\t_{n}}[I^h(\zeta_s; s,Z_s)- \eta_s] ds+ \ol u(\t_n)\1_{\{\t_n<T\}}  -\sum_{i=1}^{n-1} c^P_{\t_i} \1_{\{\t_i<T\}}\Big] \le u(\th; t, x). 
\eeaa
Send $n\to\infty$, by \reff{taun} we have
\beaa
 \dbE^{\zeta,Z}\Big[\int_t^T [ I^h(\zeta_s; s,Z_s) - \eta_s] ds  -\sum_{i=1}^\infty  c^P_{\t_i} \1_{\{\t_i<T\}}\Big] \le u(\th; t, x). 
\eeaa
Now by the arbitrariness of $(\eta, Z, \zeta)$,  we prove $u_\infty(\th; t, x)\le u(\th; t, x)$.
\qed

\ms
We conclude this section with the following result which shows that the principal agent problem \reff{VPinfty} is time consistent in the following sense.
\begin{thm}
\label{thm-DPP}
The principal agent problem \reff{VPinfty} satisfies the following dynamic programming principle: for any $\th\in \Th$ and $t<T$, and abbreviating the notation $\t^\th:= \t^{\th; t,R^\th_t,\eta, Z}$,
\bea
\label{VPDPP}
\left.\ba{c}
\dis V^P_\infty(\th; t) = \sup_{(\eta, Z)\in \cA^\th_{t,R^\th_t}}  \dbE^{ \th,Z}\Big[  \int_t^{\t^\th} [ I^h(\th; r,Z_r)-\eta_r] dr  +\sup_{\tilde \th\in \Th} V^P_\infty(\tilde\th; \t^\th) - c^P_{\t^\th}\1_{\{\t^\th<T\}}\Big];\\
\dis V^P_\infty(t) = \sup_{\th\in \Th}\sup_{(\eta, Z)\in \cA^\th_{t,R^\th_t}}  \dbE^{\th,Z}\Big[  \int_t^{\t^\th} [ I^h(\th; r,Z_r)-\eta_r] dr  + V^P_\infty(\t^\th) - c^P_{\t^\th}\1_{\{\t^\th<T\}}\Big].
\ea\right.
\eea
Moreover, let $(\eta^*, Z^*, \zeta^*)\in \cA^{\infty, \th}_{t, R^\th_t}$ with corresponding $(\{\t^*_i\}_{i\ge 0}, \cX^*)$  be an optimal control for $V^P_\infty(\th; t) = u_\infty(\th; t, R^\th_t)$, then for any $i\ge 1$, $(\eta^*, Z^*, \zeta^*)|_{[\t^*_i, T]}$ remains optimal for the problem $V^P_\infty(\zeta^*_{\t^*_i}; \t^*_i)$.
\end{thm}
\proof Recall Proposition \ref{prop-conv} (i) and \reff{oluinfty}. The first line of \reff{VPDPP} follows directly from \reff{uinftyDPP} by setting $\t\equiv T$, which implies further the second line of \reff{VPDPP}. We shall emphasize that, \reff{VPDPP} involves the term $V^P_\infty(\tilde\th; \t^{\th; t,R^\th_t,\eta, Z})$, not $V^P_\infty(\tilde\th; \t^{\tilde\th; t,R^{\tilde \th}_t,\eta, Z})$.

Moreover, if $(\eta^*, Z^*, \zeta^*)\in \cA^{\infty, \th}_{t, R^\th_t}$   is optimal for $V^P_\infty(\th; t)$, namely
\beaa
V^P_\infty(\th; t) = \dbE^{\zeta^*,Z^*}\Big[  \int_t^T [ I^h(\zeta^*_r; r,Z^*_r) -\eta^*_r] dr  -\sum_{i=1}^\infty  c^P_{\t^*_i}\1_{\{\t^*_i<T\}}\Big],
\eeaa
then $\tilde \th = \zeta^*_{\t^*_1}$ is optimal in \reff{VPDPP}, and \reff{VPDPP} becomes
\beaa
 V^P_\infty(\th; t) = \dbE^{\zeta^*,Z^*}\Big[ \int_t^{\t_1^*}[I^h(\th; r,Z^*_r) -\eta^*_r] dr  + V^P_\infty(\zeta^*_{\t^*_1}; \t^*_1) - c^P_{\t^*_1}\1_{\{\t^*_1<T\}}\Big].
\eeaa
This implies that
\beaa
 \dbE^{\zeta^*,Z^*}\Big[ V^P_\infty(\zeta^*_{\t^*_1}; \t^*_1) \Big]=\dbE^{\zeta^*,Z^*}\Big[  \int_{\t^*_1}^T[I^h(\zeta^*_r; r,Z^*_r)-\eta^*_r] dr  -\sum_{i=2}^\infty c^P_{\t^*_i}\1_{\{\t^*_i<T\}}\Big].
 \eeaa
 On the other hand, by \reff{uinfty} it is clear that, on $\{\t^*_1<T\}$,
 \beaa
V^P_\infty(\zeta^*_{\t^*_1}; \t^*_1) \ge \dbE_{\cF_{\t^*_1}}^{\zeta^*,Z^*}\Big[  \int_{\t^*_1}^T [I^h(\zeta^*_r; r,Z^*_r)-\eta^*_r] dr  -\sum_{i=2}^\infty c^P_{\t^*_i}\1_{\{\t^*_i<T\}}\Big].
 \eeaa
 Thus equality holds a.s. at above. That is, $(\eta^*, Z^*, \zeta^*)|_{[\t^*_1, T]}$ remains optimal for the problem $V^P_\infty(\zeta^*_{\t^*_1}; \t^*_1)$. Repeat the arguments we see that $(\eta^*, Z^*, \zeta^*)|_{[\t^*_i, T]}$ remains optimal for the problem $V^P_\infty(\zeta^*_{\t^*_i}; \t^*_i)$ on $\{\t^*_i<T\}$ for all $i\ge 1$.
\qed

\begin{rem}
\label{rem-consistency}
Theorem \ref{thm-DPP} shows that the principal's problem $V^P_\infty(t)$ is time consistent at the times $\t^*_i$ when the agents quit the job. We should note though, if the principal reconsiders the contract within the time period $(\t^*_i, \t^*_{i+1})$, during which the agent is staying for the job, she may find $(\eta^*, Z^*, \zeta^*)$ not optimal. However, since the principal is required to commit to the contract, she is not allowed to fire the agent and sign a better contract during those time periods, so such time inconsistency for the principal is irrelevant in our model.    
\end{rem}

\section{Comparison with the standard contracts}
\label{sect-standard0}
\setcounter{equation}{0}
In this section we investigate briefly the standard principal agent problem \reff{VP0standard} with full commitment. 
In particular, the agent is not allowed to quit, rather than having no incentive to quit as in \reff{VP0}. We shall compare the principal's optimal values under different settings.

\subsection{The dynamic value function for the standard principal agent problem}
\label{sect-standard}
We first introduce the principal's dynamic value function corresponding to \reff{VP0standard}: 
\bea
\label{u0}
\left.\ba{c}
\dis u^S_0(\th; t, x) :=  \sup_{(\eta, Z)\in \cA^\th_{t,x}} J_P(\th; t,  \eta),\q 0\le t<T, x\le \ol L^\th_t.
 \ea\right.
 \eea
We then have the following results. Since the proofs are similar, actually easier than those in Sections \ref{sect-1quit} and \ref{sect-u1}, we omit them. 
 
 \begin{thm}
 \label{thm-u0} 
 Let Assumption \ref{assum-uAh} hold. 
 
 (i) For each $\th\in \Th$ and $t< T$, $u^S_0(\th; t, \cd)$ is decreasing in $x\in (-\infty, \ol L^\th_t]$. Consequently, 
 \bea
 \label{VSP0-uS0}
\dis V^{S,P}_0(\th; t) = u^S_0(\th; t, R^\th_t),\q\mbox{provided}\q R^{\th_0}_t \le \ol L^{\th_0}_t.
 \eea

 (ii) For any $t_1, t_2 <T$ and $x_1 \le \ol L^\th_{t_1}<0$, $x_2 \le \ol L^\th_{t_2}$, we have
 \bea
 \label{uS0reg}
 \big|u^S_0(\th; t_1, x_1) - u^S_0(\th; t_2, x_2)\big| \le  C\Big[\sqrt{|x_1-x_2|} + \sqrt{\rho(|t_1-t_2|) - x_1 |t_1-t_2|}~\Big].
 \eea
 
 (iii) $u^S_0(\theta;\cd,\cd)$ has boundary conditions:
\bea
\label{uS0boundary1}
u^S_0(\theta; t, \ol L^\th_t)=-C_0(T-t),\q \forall t\in [0, T);\qq u^S_0(\theta;T-,x)=0,\q \forall x< 0.
\eea
Moreover, for $x\to-\infty$, $u_0$ satisfies the following growth condition: for $\forall t<T$, $x\le \ol L^\th_t$,
\bea
\label{uS0boundary2}
-C_0(T-t) \leq u^S_0(\theta; t,x)\leq  C\sqrt{\rho(T-t)-x(T-t)},
\eea

(iv) $u^S_0(\theta; \cd, \cd)$ is the unique viscosity solution of the following HJB equation with boundary conditions \reff{uS0boundary1} and growth conditions \reff{uS0boundary2}:
\bea
\label{uS0HJB}
\cL^\th u^S_0(\th; t,x) =0,\q 0\le t<T, x\le \ol L^\th_t.
\eea
 \end{thm}
 
Similarly to \reff{un} and \reff{VPn}, we now consider the principal agent problem where the agent is allowed to quit, however, after seeing $n$  agents' quitting, the principal will offer only committed contracts so that the new agent cannot quit anymore: for $n=1,2,\cds$, 
\bea
\label{uSn}
\left.\ba{c}
\dis \ol u^S_{n-1}(t) := \sup_{\th\in \Th} u^S_{n-1}(\th; t, R^\th_t) - c^P_t;\\
\dis u^S_n(\th; t,x) :=\!\!\!\! \sup_{(\eta, Z)\in \cA^\th_{t,x}} \!\!\!\!\dbE^{ \th,Z}\Big[\int_t^{\t^{\th; t, x, \eta, Z}} \!\!\!\!\!\!\!\!\!\! [I^h(\th; s,Z_s)- \eta_s] ds+  \ol u^S_{n-1}(\t^{\th; t, x, \eta, Z}) \1_{\{\t^{\th; t, x,\eta, Z}<T\}}\Big];\ms\\
\dis V^{S,P}_n(\th; t) := \sup_{\eta\in \cA^P_t} \dbE^{ \th, Z^{\th; t, \eta}}\Big[ \int_t^{\t^{\th; t, \eta}} [ I^h(\th; s, Z^{\th; t, \eta}_s) - \eta_s] ds  +\ol u^S_{n-1}(\t^{\th; t, \eta}) \1_{\{\t^{\th; t, \eta}<T\}}\Big],\\
\mbox{subject to}\q V^A_1(\th; t, \eta) \ge R^\th_t.
\ea\right.
\eea
Then Theorem \ref{thm-HJBn} remain true for $u^S_n, \ol u^S_n$ and $V^{S,P}_n$, after obvious modifications. Moreover, for the same function $u_\infty$ in \reff{uinfty}, following similar arguments as in Proposition \ref{prop-conv} we have

\begin{prop}
\label{prop-convS}
Under Assumptions \ref{assum-uAh} and \ref{assum-RCth}, we have $|u^S_n-u_\infty|\le {C\over n}$.
\end{prop}

\begin{rem}
\label{rem-u0}
We may also define a value function $u_0$ to characterize $V^P_0$. Denote 
 \bea
 \label{olRth}
\left.\ba{c}
\dis \ol \cA^\th_{t,x}:= \Big\{(\eta, Z)\in \cA^\th_{t,x}: \cX^{\th; t,x, \eta, Z}_s \ge R^\th_s, ~s\in [t, T], \mbox{a.s.}\Big\};\q \ol R^\th_t := \inf\Big\{x: \ol \cA^\th_{t,x} \neq \emptyset\Big\};\ms\\
\dis u_0(\th; t,x) := \sup_{(\eta, Z)\in \ol \cA^\th_{t,x}} J_P(\th; t,  \eta),\q t<T, x\in [\ol R^\th_t, \ol L^\th_t].
\ea\right.
 \eea
  Then one can easily show that $u_0 \le u^S_0$, and $u_0$ is decreasing in $x$, and thus $V^P_0(\th; t) = u_0(\th; t, \ol R^\th_t)$. However, it seems hard to establish the regularity of $u_0$ and its boundary condition at $\ol R^\th$. Nevertheless, we have the desired regularity for $V^P_0(\th; t)$ in \reff{VP0-reg}. 
\end{rem}

\subsection{Some comparisons}
\label{sect-comparison}
In this subsection we compare the principal's optimal utilities under different settings.  First, recall \reff{VP0-C0}, it is obvious that 
\beaa
V^{S,P}_0(\th; t) \ge V^P_0(\th; t),\q\mbox{and hence}\q u^S_n(\th; t,x) \ge u_n(\th; t, x)~\mbox{for all}~ \th \in \Th, (t,x)\in D_\th. 
\eeaa
That is, the principal would endure some loss by providing the agent the incentives to stay in the contract.

We next consider the case that $\Th = \{\th_0\}$ is a singleton, namely there is only one type of agent in the market, as in most publications in the literature. Recall Proposition \ref{prop-unmon}.

\begin{prop}
\label{prop-uSnmon}
Let  Assumptions \ref{assum-uAh} and \ref{assum-RCth} hold and $\Th = \{\th_0\}$. Then $u^S_n \downarrow u_\infty$ and   $V^{S,P}_n\downarrow V^P_\infty$, as $n\to \infty$.
\end{prop}
We again postpone the proof to Appendix.

\begin{rem}
\label{rem-uSnmon}
By Propositions \ref{prop-unmon} and \ref{prop-uSnmon} we see that, in the case $\Th=\{\th_0\}$, it is most desirable for the principal to offer committed contracts, where the agent is not allowed to quit, and receives the optimal utility $V^{S,P}_0$. However, if the principal cannot force the agent to commit, she would still prefer to offer a standard contract than a self-enforcing one, being aware that she might need to hire a new agent if the current one quits the job (or improve the contract by a retention offer when the agent threatens to quit).

We should note though that, in practice, if the agent is asked to commit to the contract for the whole period, he might ask for a larger  individual reservation $R^\th$. This will of course change our analysis.
\end{rem}

\ms
We now turn to the case that there are multiple (or even infinitely many) types of agents in the market, which is one of the main features of our model and is often the case in practice.  In this case, the principal may have the opportunity to choose a cheaper agent when the current one quits the job. If this gain is larger than the costs $c^{\th_0}, c^P$ induced by the agent's quitting, namely (again this is impossible when $\Th = \{\th_0\}$ as we will see in \reff{u0quit} below in the proof of Proposition \ref{prop-uSnmon}):
\beaa
\ol u_0(\t_1) = \sup_{\th\in \Th} u_0(\th;\t_1, R^\th_{\t_1}) - c^P_{\t_1} > u_0(\th_0; \t_1, R^{\th_0}_{\t_1} - c^{\th_0}_{\t_1}),
\eeaa
then  the principal may design the contract to induce the agent to quit. The following example shows that this can indeed be the case.

\begin{eg}
\label{eg-comparison}
Set $T=1$, $\Th=\{\th_0,\th_1\}\subset (0,\infty)$, $f(\th; t,\eta)= - e^{-\eta}$, $h(\th; t,\a)={1\over 2\theta}\a_t^2$,  and $R^{\th_0}_0=R^{\th_1}_0 = x_0:= -2e^{-C_0}$, where $C_0>0$ is the fixed constant for $\cA^P$. These determine the value $u^S_0(\th; 0, x_0)$ and assume without loss of generality that $u^S_0(\th_0; 0, x_0) \ge u^S_0(\th_1; 0, x_0)$. Set further that $c^\th\equiv e^{-C_0}$, $c^P\equiv 0$, and, for some constant $n\ge 1$,
\beaa
&\dis R^{\th_0}_t := [x_0 + 3e^{-C_0}t] \1_{[0, {1\over 2})}(t) - e^{-C_0}[1-t]\1_{[{1\over 2}, 1]}(t);\\
&\dis   R^{\th_1}_t := [x_0 - n t]\1_{[0, {1\over 2})}(t) + 2[ x_0 - {n\over 2}](1-t)\1_{[{1\over 2}, 1]}(t).
\eeaa
Then, $\lim_{n\to \infty} u^S_1(\th_0; 0, x_0) = \infty$, and consequently, when $n$ is large enough, we have
\bea
\label{comparison}
\max_{\th \in \Th} V^{S,P}_0(\th; 0) = u^S_0(\th_0; 0, x_0) < u^S_1(\th_0; 0, x_0) \le \max_{\th \in \Th} V^{S,P}_1(\th; 0).
\eea
\end{eg}

\proof First, it is clear that $R^{\th_1}_t \le R^{\th_0}_t \le -e^{-C_0}(1-t)= \ol L^{\th_0}_t=\ol L^{\th_1}_t$. Note that 
\beaa
\ul L^{\th_0}_t = [-2e^{-C_0} + 3e^{-C_0}t] \1_{[0, {1\over 2})}(t) - e^{-C_0}[1-t]\1_{[{1\over 2}, 1]}(t) - e^{-C_0}.
\eeaa
 Set
\beaa
\eta_t := C_0\1_{[0, {1\over 2})}(t) + [C_0-\ln 3]\1_{[{1\over 2}, 1]}(t),\q \tilde \eta \equiv -\ln (n+ 4e^{-C_0}). 
\eeaa
Then one can easily see that $(\eta, 0)\in \cA^{\th_0}_{0, x_0}$ with $\t^{\th_0; 0, x_0, \eta, 0} = {1\over 2}$, and $(\tilde \eta, 0) \in \cA^{\th_1}_{{1\over 2}, x_0 - {n\over 2}}$. Thus
\beaa
\dis u_1(\th_0; 0, x_0) &\ge& J_1(\th_0; 0, x_0, \eta, Z) = \int_0^{1\over 2} [\th_0 Z_s -\eta_s] ds +  \ol u_0({1\over 2})\\
&\ge& -{1\over 2} C_0 + u_0(\th_1; {1\over 2}, R^{\th_1}_{1\over 2}) = -{1\over 2} C_0 + u_0(\th_1; {1\over 2}, x_0 - {n\over 2});\\
u_0(\th_1; {1\over 2}, x_0 - {n\over 2}) &\ge& J_0(\th_1; {1\over 2},  \tilde \eta, \tilde Z)  = \int_{1\over 2}^1 [\ln (n+ 4e^{-C_0})] ds = {1\over 2} \ln (n+ 4e^{-C_0}).
\eeaa
 Therefore,
 \beaa
 u_1(\th_0; 0, x_0) \ge  {1\over 2} \ln (n+ 4e^{-C_0}) -{1\over 2} C_0 \to \infty,\q\mbox{as}\q n\to \infty.
 \eeaa
  
 \vspace{-0.9cm}
 \qed
 
 \ms
 \begin{rem}
\label{rem-betterquit}
In Example \ref{eg-comparison}, the agent $\th_1$ becomes cheaper when time $t$ evolves (before ${1\over 2}$). While it is desirable for the principal to hire agent $\th_0$ at initial time $0$, it becomes obvious that at later time the agent $\th_1$ is more desirable. Since in our model the principal is not allowed to fire the agent $\th_0$, she may design the contract to induce the agent $\th_0$ to quit. This will provide the principal the opportunity to hire the cheaper agent $\th_1$ at a later time and obtain a larger utility. 
\end{rem}

\section{Summary}
\label{sect-summary}
\setcounter{equation}{0}
{ In this paper we study a principal agent problem with moral hazard in finite time horizon, where only the principal is committed to the contract, but the agent is allowed to quit. Instead of offering self-enforcing contracts, the principal allows the agent to quit. When the current agent quits the job,  both the principal and the agent would incur some costs, and the principal would hire a new agent, satisfying the agent's individual rationality constraint at that time. One feature of our model 
is that there are a family of agents with different types in the market. However, unlike the third best case with adverse selection, here the principal knows the type of the agents and at any time the principal hires only one agent.

The principal's optimal utility is a function of the agent's continuation utility, which is continuous in the interior of the domain, but could be discontinuous at the boundary. Mathematically our main goal is to characterize the principal's dynamic optimal utility function through a system of HJB equations, parametrized by the agent's type.  Our results verify that the self-enforcing contracts can only be sub-optimal for the principal. That is, it is not desirable or say too expensive for the principal to disincentivize the agent from quitting. 
 It is somewhat surprising though that, in some markets the standard optimal contract for committed agents may also be sub-optimal. That is, the principal may prefer the agent to quit so that she can hire a cheaper agent from the market. In this case she would design the contract to induce the current agent to quit. This is mainly due to the presence of various types of agents. Moreover, by assuming uniform lower bound of the agent's quitting cost, the agents won't quit too frequently and  the principal will only see finitely many quittings, which is consistent with common practice. We would also like to mention that the principal's problem is time consistent in the sense that, the optimal contract remains optimal at any optimal quitting time of the agents.

We should note that, for technical reasons, in this paper each agent's dynamic individual rationality is given exogenously and is deterministic. It will be very interesting to provide endogenous models for this individual rationality, through the agent's past performance and/or through multiple principals' competition. Another interesting problem is that the principal can also fire the agent, and thus it becomes a Dynkin game problem, mixed with stochastic controls. We shall leave these important issues to  future research. 
}

\section{Appendix}
\label{sect-appendix}
\setcounter{equation}{0}

In this Appendix we present the postponed technical proofs.

\ms
\no{\bf Proof of Lemma \ref{lem-YZetaest}.}  Denote $(Y, Z) := (Y^{\th, \eta}, Z^{\th, \eta})$ for notational simplicity. Note that 
\beaa
Y_t =  \int_t^T \big[ f(\th; s,\eta) + \tilde H(\th; s,Z_s)\big] ds - \int_t^T Z_s dB^{\th,Z}_s.
\eeaa
Then, by \reff{Ltht},
 \beaa
 \ol L^\th_t - Y_t =  \dbE^{\th,Z}_t\Big[\int_t^T \big[[f(\th; s,C_0) -f(\th; s,\eta_s)]-\tilde{H}(\th; s,Z_s)\big]ds\Big].
 \eeaa
 Since $\eta\le C_0$ and $\pa_{\eta} f \ge {1\over \L_0}$, and by \reff{HIhproperty},   we have
\beaa
f(\th; s,C_0)-f(\th; s,\eta_s)\geq {1\over \L_0}(C_0-\eta_s)\ge 0;\qq -\tilde{H}(\th; s,Z_s) \ge {1\over C} |Z_s|^2 \ge 0.
\eeaa
 Then
\beaa
&&\dis \dbE^{ \th,Z}_t\Big[\int_t^T |Z_s|^2 ds\Big]\le C \dbE^{ \th,Z}_t\Big[\int_t^T[-\tilde{H}(\th; s,Z_s)] ds\Big] \le C[\ol L^\th_t - Y_t];\\
&&\dis \dbE^{\th,Z}_t\Big[\int_t^T[C_0-\eta_s] ds\Big]\le C \dbE^{ \th,Z}_t\Big[\int_t^T \big[f(\th; s,C_0) -f(\th; s,\eta_s)\big] ds\Big] \le C[\ol L^\th_t - Y_t].
\eeaa

Moreover, by \reff{HIhproperty} we have $|\eta_s|^2 \le C[1-f(\th; s, \eta_s)]$. Then, 
\beaa
&&\dis \dbE^{ \th,Z}_t\Big[\int_t^T |\eta_s|^2 ds\Big] \le C\dbE^{ \th,Z}_t\Big[\int_t^T [1-f(\th; s, \eta_s)] ds\Big] \\
&&\dis \le C\dbE^{ \th,Z}_t\Big[\int_t^T [1-f(\th; s, \eta_s) - \tilde H(\th; s, Z_s)] ds\Big] = C[T-t -Y_t].
\eeaa

\vspace{-9mm}
\qed

\bs\bs

\no {\bf Proof of Propostion \ref{prop-VP0}.} Without loss of generality we estimate only $|V^P_0(\th; 0) - V^P_0(\th; \d)|$ for $\d\le T$.

{\it Step 1.} Let $\eta\in \cA^\dbP=\cA^P_0$ satisfies the constraint in \reff{VP0}, and simplify the notation $(Y, Z) := (Y^{\th, \eta}, Z^{\th, \eta})$. Clearly $\eta_{[\d, T]}$ also satisfies the constraint in \reff{VPt} with $t=\d$, then 
\beaa
\dbE^{\th, Z}_{\cF_\d}\Big[  \int_\d^T \big[I^h(\th; s,Z_s)-\eta_s\big] ds\Big] \le V^P_0(\th; \d),\q\mbox{a.s.}
\eeaa
Thus, by \reff{HIhproperty} and \reff{YZetaest},
\beaa
&&\dis J_P(\th; 0, \eta) - V^P_0(\th; \d) \le \dbE^{\th, Z}\Big[  \int_0^\d \big[I^h(\th; s,Z_s)-\eta_s\big] ds\Big]\\
&&\dis \le \dbE^{\th, Z}\Big[  \int_0^\d \big[C|Z_s|+|\eta_s|\big] ds\Big] \le C\sqrt{\d} \Big(\dbE^{\th, Z}\Big[  \int_0^T \big[|Z_s|^2+|\eta_s|^2\big] ds\Big]\Big)^{1\over 2}\\
&&\dis \le C\sqrt{\d} \Big( \ol L^\th_0 - Y_0 + T-Y_0 \Big)^{1\over 2}\le   C\sqrt{\d} \Big( \ol L^\th_0  + T-2R^\th_0 \Big)^{1\over 2} \le C\sqrt{\d}.
\eeaa
Since $\eta$ is arbitrary, we obtain $V^P_0(\th; 0)-V^P_0(\th; \d)\le C\sqrt{\d}$.

{\it Step 2.} Let $\eta\in \cA^\dbP_\d$ satisfies the constraint in \reff{VPt}. If suffices to prove
\bea
\label{JP<V0}
J_P(\th; \d, \eta) - V^P_0(\th; 0)  \le C\sqrt{\rho(\d)}.
\eea
Then by the arbitrariness of $\eta$ we obtain $V^P_0(\th; \d)-V^P_0(\th; 0)\le C\sqrt{\rho(\d)}$. 

Denote $\d':=  \rho(\d)+ \L_0\d$ and  $(Y, Z) := (Y^{\th, \eta}, Z^{\th, \eta})$. We prove \reff{JP<V0} in two cases.

{\it Case 1.} $Y_\d \ge \ol L^\th_\d - \d' $. Then, by \reff{HIhproperty}, \reff{YZetaest}, and \reff{VP0-C0},
\beaa
J_P(\th; \d, \eta) &=& \dbE^{\th, Z}\Big[  \int_\d^T \big[I^h(\th; s,Z_s)-\eta_s\big] ds\Big] \\
&\le& \dbE^{\th, Z}\Big[  \int_\d^T \big[C|Z_s| + C_0-\eta_s\big] ds\Big] - C_0(T-\d)\\
&\le& C\Big[\sqrt{\ol L^\th_\d - Y_\d} + \ol L^\th_\d - Y_\d\Big] + V^P_0(\th; 0) + C_0\d \\
&\le& C\Big[\sqrt{\d'} +\d'\Big] + V^P_0(\th; 0) + C_0\d ~\le~  V^P_0(\th; 0)  + C\sqrt{\rho(\d)}.
\eeaa

{\it Case 2.} $Y_\d < \ol L^\th_\d -\d'$. Construct 
\beaa
&\dis \t_\d:= \inf\big\{s\ge \d: Y_s = \ol L^\th_s - \d'\big\}<T,~ s\in [\d, T],\\
&\dis \eta^\d_s:= C_0 \1_{[0, \d)}(s) + \eta_s \1_{[\d, \t_\d)}(s) + C_0\1_{[\t_\d, T]}(s).
\eeaa
It is clear that $\eta^\d\in \cA^P$. Denote $(Y^\d, Z^\d):= (Y^{\th, \eta^\d}, Z^{\th, \eta^\d})$. One can easily verify that
\beaa
(Y^\d_s, Z^\d_s) = (Y_s + \d', Z_s) \1_{[\d, \t_\d]}(s) + (\ol L^\th_s, 0)\1_{(\t_\d, T]}(s),\q s\in [\d, T].
\eeaa
In particular, $Y^\d_s \ge Y_s\ge R^\th_s$ for $s\in [t, T]$. Moreover, for $s\in [0, \d]$, we have
\beaa
Y^\d_s = Y^\d_\d  + \int_t^\d \big[ f(\th; s,\eta^\d) + H(\th; s,Z^\d_s)\big] ds - \int_t^\d Z^\d_s dB_s.
\eeaa
Since $\eta\in \cA^P_\d$ is $\dbF^{B^\d}$-progressively measurable,  then $Y_\d$ is deterministic, and hence so is $Y^\d_\d$. This implies $Z^\d_s =0$ for $s\in [0, \d]$. Therefore, by \reff{uA} and the regularity of $R^\th$,
\beaa
Y^\d_s = Y^\d_\d  + \int_s^\d f(\th; r, C_0) dr \ge  Y^\d_\d  - \L_0\d  = Y_\d + \rho(\d) \ge R^\th_\d + \rho(\d) \ge R^\th_s,\q s\in [0, \d].
\eeaa
That is, $\eta^\d$ satisfies the constraints in \reff{VP0}. Then
\beaa
&&\dis J_P(\th; \d, \eta) - V^P_0(\th; 0) \le J_P(\th; \d, \eta) - J_P(\th; 0, \eta^\d) \\
&&\dis = \dbE^{\th, Z}\Big[  \int_\d^T \big[I^h(\th; s,Z_s)-\eta_s\big] ds\Big] -  \dbE^{\th, Z^\d}\Big[ \int_\d^{\t_\d} \big[I^h(\th; s,Z_s)-\eta_s\big] ds + C_0[T-\t_\d-\d]\Big].
\eeaa
Since $Z^\d = Z$ on $[\d, \t_\d]$, then $\dbP^{\th, Z^\d} = \dbP^Z$ on $\cF^\d_{\t_\d}$. Thus, by \reff{YZetaest},
\beaa
&&\dis J_P(\th; \d, \eta) - V^P_0(\th; 0)  \\
&&\dis \le \dbE^{\th, Z}\Big[  \int_\d^T \big[I^h(\th; s,Z_s)-\eta_s\big] ds\Big] -  \dbE^{\th, Z}\Big[ \int_\d^{\t_\d} \big[I^h(\th; s,Z_s)-\eta_s\big] ds + C_0[T-\t_\d-\d]\Big]\\
&&\dis \le \dbE^{\th, Z}\Big[  \int_{\t_\d}^T \big[I^h(\th; s,Z_s)+C_0-\eta_s\big] ds\Big] \le \dbE^{\th, Z}\Big[  \int_{\t_\d}^T \big[C|Z_s|+C_0-\eta_s\big] ds\Big]\\
&&\dis  \le C\dbE^{\th, Z}\Big[ \sqrt{\ol L^\th_{\t_\d} - Y_{\t_\d}} + \ol L^\th_{\t_\d} - Y_{\t_\d} \Big]=C\Big[\sqrt{\d'} +  \d' \Big]\le C\sqrt{\rho(\d)}.
\eeaa

\vspace{-9mm}
\qed

\bs\bs

To prove Proposition  \ref{prop-agent1}, we shall first prove Lemma 3.4.

\no {\bf Proof of Lemma \ref{lem-u1domain}.} (i) Recall \reff{etan} and fix $n\ge 1$ sufficiently large.  Clearly $\eta^n$ is bounded from below, and  $\pa_yI^f> 0$ implies $\eta^n_s \le  C_0$, thus $\eta^n\in \cA^P_t$. Denote 
\beaa
\cX^n_s := x - \int_t^s f(\th; r, \eta^n_r) dr = x - \int_t^s \big[f(\th; r, C_0) - n\big]dr = x - \ol L^\th_t + \ol L^\th_s + n(s-t),
\eeaa
for $s\in [t, t+\e_n]$. Then it is clear that $\cX^n_s < \ol L^\th_s$ for $s\in [t, t+\e_n)$ and $\cX^n_{t+\e_n} = \ol L^\th_{t+\e_n}$. Since $\eta^n_s = C_0$ for $s> t+\e_n$, this implies further that $\cX^n_s = \ol L^\th_s$ for $s\in [t+\e_n, T]$. In particular, $\cX^n_T = \ol L^\th_T=0$, and thus $(\eta^n, 0) \in \cA^\th_{t,x}$. 

Moreover, for $s\in [t, t+\e_n]$, by \reff{HIhproperty} and \reff{uA} we have
\beaa
|\eta^n_s|^2 \le C[1- f(\th; s, \eta^n_s)] = C [1 - f(\th; s, C_0) + n]\le C [1 +\L_0 + n] \le Cn.
\eeaa
 Then
\beaa
\int_t^T [C_0 - \eta^n_s]ds= \int_t^{t+\e_n} [C_0 - \eta^n_s]ds \le \int_t^{t+\e_n}[C_0+ |\eta^n_s|]ds \le [{C_0\over n} + {C\over \sqrt{n}}] (\ol L^\th_t - x).
\eeaa
This verifies  \reff{u1domain2} immediately by choosing $n$ large enough.

(ii) For $s\in [t, t+\e_n]$, by \reff{uA} and  Assumption \ref{assum-RCth} (iii) we have, whenever $n\ge \L_1$,
\beaa
\cX^n_s - \ul L^\th_s  > \ul L^\th_t  - \int_t^s \big[f(\th; r, C_0) - n\big]dr- \ul L^\th_s \ge  n(s-t)- |\ul L^\th_s -\ul L^\th_t| \ge 0.
\eeaa 
Moreover, for $s\in (t+\e_n, T]$, we have $\cX^n_s = \ol L^\th_s > \ul L^\th_s$.   So $(\cX^n, 0, 0)$ satisfies RBSDE \reff{RBSDE} with $\eta^n$, and thus  $\t^{\th; t, \eta^{\th;t,x}}=T$.
\qed

\bs

\no{\bf Proof of Proposition \ref{prop-agent1}.} (i) First,  for fixed $\t\in \cT_t$, the agent's optimal effort can be solved through a BSDE:
\bea
\label{BSDE1}
Y^{\th; t,\eta, \t}_s = \ul L^{\th}_\t \1_{\{\t<T\}} + \int_s^\t \big[ f(\th; r,\eta_r) + H(\th; r,Z^{\th;t,\eta,\t}_r)\big] dr - \int_s^\t Z^{\th; t,\eta,\t}_r dB_r.
\eea
Similarly to Proposition \ref{prop-agent}, we see that $V^A_1(\th; t, \eta) = \sup_{\t\in \cT_t} Y^{\th; t,\eta, \t}_t$ and, once $\t$ is given, the agent's optimal effort is $\a^* =I^h(\th; s,Z^{\th;t,\eta,\tau})$ on $[t, \t]$. Then it follows from the standard RBSDE theory that $V^A_1(\th; t, \eta) = Y^{\th; t, \eta}_t$ and $\t^{\th; t,\eta}$ is the smallest optimal stopping time. Notice further that $K^{\th; t, \eta}_s =0$ for $t\le s\le \t^{\th;t, \eta}$, then \reff{RBSDE} coincide with \reff{BSDE1}  on $[t, \t^{\th; t, \eta}]$. Thus $Z^{\th; t, \eta}_s = Z^{\th; t,\eta,\t^{\th; t, \eta}}_s$, $s\in [t, \t^{\th; t, \eta}]$. Since $\ul L^{\th}_\t \1_{\{\t<T\}}$ is bounded, by \reff{BSDE1} we see that $Z^{\th; t,\eta,\t}\in \cA^A_t$ (see e.g. \cite[Chapter 7]{Zhang}), and hence $\a^{\th; t,\eta} = I^h(\th; s,Z^{\th;t,\eta})\1_{[t,\tau^{\th;t,\eta}]}\in \cA^A_t$.

(ii) First, by definition $Y^{\th; t, \eta}_t \ge \ul L^\th_t$. Next, since $\ol L^\th > \ul L^\th$, then $(\ol L^\th, 0, 0)$ is the unique solution of the RBSDE \reff{RBSDE} with $\eta \equiv C_0$, and thus it follows from the comparison principle of RBSDEs that $Y^{\th; t, \eta}_t \le \ol L^\th_t$. On the other hand, for any $x\in (\ul L^\th_t, \ol L^\th_t)$, by  Lemma \ref{lem-u1domain}  we see that $\eta^{\th; t,x}\in \cA^P_t$ and $x= Y^{\th; t, \eta^{\th; t,x}}$. Moreover,  $x=\ul L^\th_t$, from the arguments in Lemma \ref{lem-u1domain}  we can easily see that we still have $x= Y^{\th; t, \eta^{\th; t,x}}_t$, however, in this case $Y^{\th; t, \eta^{\th;t,x}}_s > \ul L^\th_s$ holds only for $s\in (t, T]$. Then we obtain the claimed equality. 
\qed

\bs

\no{\bf Proof of Lemma \ref{lem-olu0}.}  First, by Proposition \ref{prop-VP0} and Assumption \ref{assum-RCth}, clearly $ \ol u_0(t)$ is uniformly continuous in $[0, T)$ and hence on $[0, T]$. Next, for any $t<T$ and $(\eta, Z)\in \cA^\th_{t, R^\th_t}$,  by \reff{YZetaest} we have
\beaa
 |J_P(\th; t, \eta)| &\le& \dbE^{\th, Z}\Big[\int_t^T \big[C|Z_s| + |\eta_s|\big]ds\Big] \le C\sqrt{T-t}\Big(\dbE^{\th, Z}\Big[\int_t^T \big[|Z_s|^2 + |\eta_s|^2\big]ds\Big]\Big)^{1\over 2}\\
 &\le& C\sqrt{T-t}\Big( \ol L^\th_t - R^\th_t  + T-t - R^\th_t\Big)^{1\over 2}.
\eeaa
Then
\beaa
|V^P_0(\th; t)|\le  C\sqrt{T-t}\Big( \ol L^\th_t - R^\th_t  + T-t - R^\th_t\Big)^{1\over 2}
\eeaa
This implies $\ol u_0$ is bounded and $\lim_{t\uparrow T} |\ol u_0(t)| =0$.  
\qed

\bs

\no {\bf Proof of Lemma \ref{lem-Zetaest}.} Denote
\beaa
(\ul \eta, \ul Z):= (\ul \eta^{\th; t, x}, \ul Z^{\th; t, x}),\q \cX:= \cX^{\th; t,x , C_0, 1},\q \t:= \t^{\th; t,x, C_0, 1},\q \ul\t:= \t^{\th; t, x, \ul \eta^{\th; t, x}, \ul Z^{\th; t, x}}.
\eeaa
By Lemma \ref{lem-u1domain} one can easily verify that $(\ul \eta, \ul Z)\in \cA^\th_{t,x}$ and $\ul \t = \t \1_{\{\t \le T_\d\}} + T \1_{\{\t  > T_\d\}}$.
Then
\beaa
&&\dis\ol u_0(t) - J_1(\th; t, x, \ul \eta, \ul Z)  = \ol u_0(t) - \dbE^{\th, 1}\Big[  \ul u_0(\t) \1_{\{\t\le T_\d\}} \\
&&\dis +\1_{\{\t\le T_\d\}} \int_t^\t \big[I^h(\th; s, 1) - C_0\big] ds +  \1_{\{\t> T_\d\}}\big[ \int_t^{T_\d} \big[I^h(\th; s, 1) - C_0\big] ds - \int_{T_\d}^T \eta_s^{\th;  T_\d, \cX_{T_\d}}ds\big]\Big].
\eeaa
Note that $t\le T_{2\d}$, then $\{\t>T_\d\}\subset \{\t> t+ \d\}$. Then, by considering the two cases $\{\t\le t+ \d\}$ and   $\{\t > t+\d\}$, it follows from \reff{u1domain2} and Lemma \ref{lem-olu0} that
\beaa
\left.\ba{lll}
\dis\ol u_0(t) - J_1(\th; t, x, \ul \eta, \ul Z)  \le \rho_0(\d) +\d+ C\dbP^{\th,1}(\t>t+\d) \ms\\
\dis \le C\rho_0(\d) +   \dbP^{\th,1}\Big(\inf_{t\le s \le t+\d} [\cX_s - \ul L^\th_s] > 0\Big)  \\
 \dis =  C\rho_0(\d) +\dbP^{\th,1}\Big(\inf_{t\le s \le t+\d} \big[x - \int_t^s \big[f(\th; r,C_0)+\tilde{H}(\th; r,1) \big]dr + [B_s^{\th,1}-B_t^{\th,1}] - \ul L^\th_s\big] > 0\Big).
 \ea\right.
  \eeaa
  Note that
  \beaa
  - \int_t^s \big[f(\th; r,C_0)+\tilde{H}(\th; r,1) \big]dr  \le C \d,\q x - \ul L^\th_s = x- \ul L^\th_t + \ul L^\th_t - \ul L^\th_s \le \d + \L_1\d.
  \eeaa
  Thus
\beaa
\ol u_0(t) - J_1(\th; t, x, \ul \eta, \ul Z)  &\le& C\rho_0(\d) + \dbP^{\th,1}\Big(  \inf_{t\le s \le t+\d}  [B_s^{\th,1}-B_t^{\th,1}] > - C \d \Big) \\
&\le& C\rho_0(\d) + C\sqrt{\d} \le C\rho_0(\d).
  \eeaa

\vspace{-8mm}
\qed

\bs
\no{\bf Proof of Proposition \ref{prop-unalternative}.} The result is trivial when $n=1$. We shall prove \reff{un2} only for $n=2$. The general case $n>2$ can be proved following similar but slightly more involved arguments. Let $\tilde u_2(\th; t, x)$ denote the right side of the first line in \reff{un2} with $n=2$. 

{\it Step 1.} We first prove $\tilde u_2(\th; t, x) \le u_2(\th; t, x)$. For any $(\eta, Z, \zeta)$ satisfying the required properties in \reff{un2}, and for any $\d>0$ small, denote 
\beaa
&\dis (\eta^0_s, Z^0_s) := (\eta_s, Z_s) \1_{[t, \t_1\wedge T_\d)}(s) + (\eta^{\th; \t_1\wedge T_\d, \cX^{\th; t, x, \eta, Z}_{\t_1\wedge T_\d}}_s, ~ 0)   \1_{[\t_1\wedge T_\d, T)}(s);\\
&\dis (\eta^1_s, Z^1_s) := (\eta_s, Z_s) \1_{[\t_{1}, \t_2\wedge T_\d)}(s) + (\eta^{\zeta_{\t_1}; \t_2\wedge T_\d, \cX^{\zeta_{\t_1}; \t_1, R^{\zeta_{\t_1}}_{\t_1}, \eta, Z}_{\t_2\wedge T_\d}}_s, ~ 0)   \1_{[\t_2\wedge T_\d, T)}(s),~\mbox{on}~ \{\t_1\le T_\d\}.
\eeaa
Then 
\beaa
&\dis (\eta^0, Z^0)\in \cA^{\th}_{t, x}, \q \t^\d_1:= \t^{\th; t, x, \eta^0, Z^0} = \t_1 \1_{\{\t_1 \le T_\d\}} + T \1_{\{\t_1> T_\d\}};\\
&\dis (\eta^1, Z^1)\in \cA^{\zeta_{\t_1}}_{\t_1, R^{\zeta_{\t_1}}_{\t_1}},\q \t^\d_2:= \t^{\zeta_{\t_1}; \t_1, R^{\zeta_{\t_1}}_{\t_1}, \eta^1, Z^1} = \t_2 \1_{\{\t_2 \le T_\d\}} + T \1_{\{\t_2> T_\d\}},~\mbox{on}~ \{\t_1\le T_\d\}.
\eeaa
Thus, by \reff{un},
\beaa
u_2(\th; t, x) &\ge&\dbE^{ \th,Z^0}\Big[\int_t^{\t^\d_1} [I^h(\th; s,Z^0_s)- \eta^0_s] ds+  \ol u_1(\t_1) \1_{\{\t_1 \le T_\d\}}\Big]\\
&\ge& \dbE^{ \th,Z^0}\Big[\int_t^{\t^\d_1} [I^h(\th; s,Z^0_s)- \eta^0_s] ds+  \big[u_1(\zeta_{\t_1}; \t_1, R^{\zeta_{\t_1}}_{\t_1}) - c^P_{\t_1}\big]\1_{\{\t_1 \le T_\d\}}\Big]\\
&\ge& \dbE^{ \th,Z^0}\Big[\int_t^{\t^\d_1} [I^h(\th; s,Z^0_s)- \eta^0_s] ds\\
&&+  \Big( \dbE^{ \zeta_{\t_1},Z^1}_{\t_1}\big[\int_{\t_1}^{\t^\d_2} [I^h(\th; s,Z^1_s)- \eta^1_s] ds + \ol u_0(\t_2)\1_{\{\t_2\le T_\d\}}\big]- c^P_{\t_1}\Big)\1_{\{\t_1 \le T_\d\}}\Big]\\
&=&\dbE^{ \zeta,Z}\Big[\1_{\{\t_2\le T_\d\}}\big[\int_t^{\t_2} [I^h(\zeta_s; s,Z_s)- \eta_s] ds - c^P_{\t_1} + \ol u_0(\t_2)\big] \\
&&\dis + \1_{\{\t_1\le T_\d<\t_2\}}\big[\int_t^{T_\d} [I^h(\zeta_s; s,Z_s)- \eta_s] ds - \int_{T_\d}^T\eta^{\zeta_{\t_1}; T_\d, \cX^{\zeta_{\t_1}; \t_1, R^{\zeta_{\t_1}}_{\t_1}, \eta, Z}_{T_\d}}_s ds - c^P_{\t_1}\big]\\
&&\dis  + \1_{\{\t_1>T_\d\}}\big[\int_t^{T_\d} [I^h(\zeta_s; s,Z_s)- \eta_s] ds - \int_{T_\d}^T\eta^{\th; T_\d, \cX^{\th; t, x, \eta, Z}_{T_\d}}_s ds \big]\Big].
\eeaa
This implies that
\beaa
&&\dis J_2(t, x, \eta, Z, \zeta) - u_2(\th; t, x)\\
&&\dis \le \dbE^{ \zeta,Z}\Big[ \1_{\{\t_2>T_\d\}}\int_{T_\d}^{\t_2} [I^h(\zeta_s; s,Z_s)- \eta_s] ds  + \ol u_0(\t_2)\1_{\{T_\d<\t_2<T\}} - c^P_{\t_1}\1_{\{T_\d<\t_1<T\}} \\
&&\dis \qq\q  + \1_{\{\t_1\le T_\d<\t_2\}}\int_{T_\d}^T\eta^{\zeta_{\t_1}; T_\d, \cX^{\zeta_{\t_1}; \t_1, R^{\zeta_{\t_1}}_{\t_1}, \eta, Z}_{T_\d}}_s ds+ \1_{\{\t_1>T_\d\}} \int_{T_\d}^T\eta^{\th; T_\d, \cX^{\th; t, x, \eta, Z}_{T_\d}}_s ds\big]\Big]\\
&&\dis \le \dbE^{ \zeta,Z}\Big[ \1_{\{\t_2>T_\d\}}\int_{T_\d}^{\t_2} [I^h(\zeta_s; s,Z_s)- \eta_s] ds  + \ol u_0(\t_2)\1_{\{T_\d<\t_2<T\}} - c^P_{\t_1}\1_{\{T_\d<\t_1<T\}} \\
&&\dis \qq\q  + C_0\d \big[\1_{\{\t_1\le T_\d<\t_2\}}+ \1_{\{\t_1>T_\d\}}\big] \Big].
\eeaa
Send $\d\to 0$, we obtain $J_2(t, x, \eta, Z, \zeta) \le u_2(\th; t, x)$. Since $(\eta, Z, \zeta)$ are arbitrary, we have $\tilde u_2(\th; t, x)\le u_2(\th; t, x)$.

{\it Step 2.} We next show that $u_2(\th; t, x) \le \tilde u_2(\th; t, x)$. Let $\tilde \rho$ be a generic modulus of continuity function. Fix a small $\d>0$. First, by \reff{un} there exists $(\eta^0, Z^0)\in \cA^\th_{t,x}$ such that
\beaa
u_2(\th; t, x) \le \dbE^{\th, Z^0}\Big[\int_t^{\t_0} [I^h(\th; s,Z^0_s)- \eta^0_s] ds  + \ol u_1(\t_0)\1_{\{\t_0<T\}}\Big] + \d.
\eeaa 
where $\t_0 := \t^{\th; t,x,\eta^0, Z^0}$. Let $t=t_0<\cds<t_n=T$ be a partition of $[t, T]$ such that $t_{n-1} = T_\d$ and $t_i-t_{i-1}\le \d$, $i=1,\cds, n-1$. Now fix an $i=1,\cds, n-1$. Choose $\th_i$ such that 
\beaa
\ol u_1(t_i) \le u_1(\th_i; t_i, R^\th_{t_i}) - c^P_{t_i} + \d.
\eeaa
Denote $x_i:= \inf_{t_{i-1}\le t\le t_i} R^{\th_i}_{t_i}$. Since $R^\th_t - \ul L^\th_t=c^\th_t\ge {1\over \L_0}$, by the uniform regularity of $R^\th_t$ we have $x_i > \ul L^\th_t+ \L_0\d$ for $t_{i-1}\le t\le t_i$, provided $\d$ is small enough. Let $(\eta^i, Z^i) \in \cA^{\th_i}_{t_i, x_i}$ be such that, denoting $\t_i := \t^{\th_i; t_i, x_i, \eta^i, Z^i}$, 
\beaa
u_1(\th_i; t_i, x_i) \le  \dbE^{\th_i, Z^i}\Big[\int_{t_i}^{\t_i} [I^h(\th_i; s,Z^i_s)- \eta^i_s] ds  + \ol u_0(\t_i)\1_{\{\t_i<T\}}\Big] + \d.
\eeaa 
 Denote, on $\{t_{i-1}< \t_0\le t_i\}$, 
\beaa
x_i':= \cX^{\th_i; \t_0, R^{\th_i}_{\t_0}, C_0, 0}_{t_i} = R^{\th_i}_{\t_0} - \int_{\t_0}^{t_i} f(\th_i; s, C_0) ds .
\eeaa
We shall note that $x_i'$ is random here, and it is obvious that $0\le x_i'-x_i \le \tilde \rho(\d)$. Following the arguments in Proposition \ref{prop-u1mon} Step 1, we can construct $(\tilde \eta^i, \tilde Z^i) \in \cA^{\th_i}_{t_i, x_i'}$ such that
\beaa
&&\dis \dbE^{\th_i, Z^i}\Big[\int_{t_i}^{\t_i} [I^h(\th_i; s,Z^i_s)- \eta^i_s] ds  + \ol u_0(\t_i)\1_{\{\t_i<T\}}\Big] \\
&&\dis \le \dbE^{\th_i, \tilde Z^i}_{t_i}\Big[\int_{t_i}^{\tilde \t_i} [I^h(\th_i; s, \tilde Z^i_s)- \tilde \eta^i_s] ds  + \ol u_0(\tilde \t_i)\1_{\{\tilde \t_i<T\}}\Big] +\tilde \rho(\d),
\eeaa
where $\tilde \t_i := \t^{\th_i; t_i, x_i', \tilde \eta^i, \tilde Z^i}$, and $\tilde \eta^i$ has a lower bound. Now extend $(\tilde \eta^i, \tilde Z^i)$ to $[\t_0, T]$ with $(\tilde \eta^i, \tilde Z^i)  = (C_0, 0)$ on $[\t_0, t_i]$. Note that, for $\d>0$ small enough, $\cX^{\th_i; \t_0, R^\th_{\t_0}, C_0, 0}_s > \ul L^\th_s$, $\t_0\le s\le t_i$.  Then we have $(\tilde \eta^i, \tilde Z^i) \in \cA^{\th_i}_{\t_0, R^{\th_i}_{\t_0}}$, and $\t^{\th_i; \t_0, R^{\th_i}_{\t_0}, \tilde \eta^i, \tilde Z^i}=\tilde \t_i$, and thus we can easily check that, again on $\{t_{i-1}< \t_0\le t_i\}$, 
\beaa
&&\dis \dbE^{\th_i, Z^i}\Big[\int_{t_i}^{\t_i} [I^h(\th_i; s,Z^i_s)- \eta^i_s] ds  + \ol u_0(\t_i)\1_{\{\t_i<T\}}\Big] \\
&&\dis \le \dbE^{\th_i, \tilde Z^i}_{\t_0}\Big[\int_{\t_0}^{\tilde \t_i} [I^h(\th_i; s, \tilde Z^i_s)- \tilde \eta^i_s] ds  + \ol u_0(\tilde \t_i)\1_{\{\tilde \t_i<T\}}\Big] +\tilde \rho(\d),
\eeaa

We now define
\beaa
&\dis (\eta_s, Z_s, \zeta_s) := (\eta^0_s, Z^0_s, \th) \1_{[t, \t_0\wedge T_\d]}(s)  + \1_{[\t_0\wedge T_\d, T]}(s)\times \\
&\dis \Big[\sum_{i=1}^{n-1}(\tilde \eta^i, \tilde Z^i, \th_i) \1_{\{t_{i-1}< \t_0 \le t_i\}} + (\eta^{\th; T_\d, \cX^{\th; t,x, \eta^0, Z^0}_{T_\d}}_s, 0, \th) \1_{\{\t_0 > t_{n-1}\}}\Big].
\eeaa 
Recall Remark \ref{rem-u1domain}, and since the above expression involves only finitely many $\eta^i$, the above $\eta$ has a uniform lower bound, then one can verify that $(\eta, Z, \zeta)\in   \cA^{\infty, \th}_{t,x}$ and 
\beaa
u_2(\th; t, x) \le J_2(t, x, \eta, Z, \zeta) + \tilde \rho(\d) \le \tilde u_2(\th; t, x)+ \tilde \rho(\d).
\eeaa
Send $\d\to 0$, we obtain $u_2(\th; t, x) \le  \tilde u_2(\th; t, x)$.
\qed

\bs
\no{\bf Proof of Proposition \ref{prop-taufinite}.}  (i) Note that, for each $i\ge 1$, on $\{\t_i < T\}$ we have
\beaa
&&\dis \ul L^{\zeta_{\t_i}}_{\t_{i+1}} \1_{\{\t_{i+1}<T\}} + \cX_{T-}\1_{\{\t_{i+1}=T\}} = \cX_{\t_{i+1}-} \\
&&\dis = R^{\zeta_{\t_i}}_{\t_i} - \int_{\t_i}^{\t_{i+1}}\big[f(\zeta_{\t_i}; s,\eta_s)+\tilde{H}(\zeta_{\t_i}; s,Z_s)\big] ds + \int_{\t_i}^{\t_{i+1}} Z_s dB^{\zeta_{\t_i}, Z}_s\\
&&\dis \ge \ul L^{\zeta_{\t_i}}_{\t_i} + {1\over \L_0} + \int_{\t_i}^{\t_{i+1}} Z_r dB^{\zeta_{\t_i}, Z}_r.
\eeaa
Then
\beaa
\dbE_{\t_i}^{\zeta,Z} \Big[ \ul L^{\zeta_{\t_i}}_{\t_{i+1}} \1_{\{\t_{i+1}<T\}} + \cX_{T-}\1_{\{\t_{i+1}=T> \t_i\}}\Big] \ge \dbE_{\t_i}^{\zeta,Z} \Big[  \big[\ul L^{\zeta_{\t_i}}_{\t_i} + {1\over \L_0}\big]\1_{\{\t_i<T\}}\Big],
\eeaa
and thus, noting that $\ul L$ and $\cX$ are uniformly bounded, 
\beaa
{1\over \L_0} \dbP^{\zeta,Z}(\t_i<T) &\le& \dbE^{\zeta,Z} \Big[ \ul L^{\zeta_{\t_i}}_{\t_{i+1}} \1_{\{\t_{i+1}<T\}} - \ul L^{\zeta_{\t_i}}_{\t_i} \1_{\{\t_i<T\}} + \cX_{T-}\1_{\{\t_{i+1}=T> \t_i\}}\Big] \\
&\le&\dbE^{\zeta,Z} \Big[ \L_1 [\t_{i+1}-\t_i] \1_{\{\t_i < T\}}  + [\cX_{T-} - \ul L^{\zeta_{\t_i}}_T]\1_{\{\t_{i+1}=T> \t_i\}}\Big]\\
&\le& C \dbE^{\zeta,Z}\Big[\t_{i+1}-\t_i  + \1_{\{\t_{i+1}=T> \t_i\}}\Big].
\eeaa
Take summation over $i$, we have
\beaa
{n\over \L_0} \dbP^{\zeta,Z}(\t_n < T) &\le& {1\over \L_0} \sum_{i=1}^n \dbP^{\zeta,Z}(\t_i < T) \le C\sum_{i=1}^n \dbE^{\zeta,Z}\Big[\t_{i+1}-\t_i  + \1_{\{\t_{i+1}=T> \t_i\}}\Big]\le C.
\eeaa
This implies \reff{taun} immediately.

(ii) {\it Step 1.} First, for any $(\eta, Z, \zeta)\in \cA^{\infty, \th}_{t,x}$ as in \reff{un2} and $\d>0$ small, construct
\bea
\label{tildeinfty}
\t^\d_n:= \t_n \wedge T_\d,\q (\tilde \eta, \tilde Z, \tilde \zeta)_s :=  (\eta, Z, \zeta)_s \1_{[0, \t^\d_n)}(s) + (\eta^{\zeta_{\t^\d_n}; \t^\d_n, R^{\zeta_{\t^\d_n}}_{\t^\d_n}}_s, 0, \zeta_{\t_n^\d}) \1_{[\t_n^\d, T)}(s).
\eea
Then clearly $(\tilde \eta, \tilde Z, \tilde \zeta)\in \cA^{\infty, \th}_{t,x}$ satisfies the constraints in \reff{uinfty}, with corresponding $\tilde \t_i = \t_i$ on $\{\t_i \le T_\d\}$ for $i\le n$, and $\tilde \t_i = T$ otherwise. Then 
\beaa
&&\dis J_n(t, x, \eta, Z, \zeta) - u_\infty(\th; t,z) \le J_n(t, x, \eta, Z, \zeta) - J_\infty(t, x, \tilde \eta, \tilde Z, \tilde \zeta)\\
&&\dis = \dbE^{\zeta,Z}\Big[\int_{\t_n^\d}^{\t_n} [I^h(\zeta_s;s,Z_s) - \eta_s]ds  +\int_{\t_n^\d}^T \eta^{\zeta_{\t^\d_n}; \t^\d_n, R^{\zeta_{\t^\d_n}}_{\t^\d_n}}_s ds\\
&&\dis\qq - \sum_{i=1}^n c^P_{\t_i} \1_{\{T_\d< \t_i <T\}} + \sup_\th V^P_0(\th; \t_n)\1_{\{\t_n<T\}}\Big]\\
&&\dis \le \dbE^{\zeta,Z}\Big[\int_{\t_n^\d}^{\t_n}  [I^h(\zeta_s;s,Z_s) - \eta_s] ds  +C_0(T-\t_n^\d) + C\1_{\{\t_n<T\}}\Big].
\eeaa
Send $\d\to 0$, one can easily see that
\beaa
 J_n(t, x, \eta, Z, \zeta) - u_\infty(\th; t,z) \le  \dbE^{\zeta,Z}\Big[C_0(T-\t_n) + C\1_{\{\t_n<T\}}\Big] \le C\dbP^{\zeta,Z}(\t_n<T) \le {C\over n},
\eeaa
 and thus
$
 u_n(\th; t, x) - u_\infty(\th; t,z) \le  {C\over n}.
$
 
 {\it Step 2.} On the other hand, for any $(\eta, Z, \zeta)\in \cA^{\infty, \th}_{t,x}$ satisfying the constraints in \reff{uinfty}, again let $(\tilde \eta, \tilde Z, \tilde \zeta)$ be constructed by \reff{tildeinfty}. Then $(\tilde \eta, \tilde Z, \tilde \zeta)\in \cA^{\infty, \th}_{t,x}$ satisfies the constraints in \reff{un2} with $\tilde \t_i = \t_i$ on $\{\t_i \le T_\d\}$,  $i=0,\cds, n$. Thus, by \reff{u1domain2},
 \beaa
&&\dis  J_\infty(t, x,  \eta,  Z,  \zeta) - u_n(\th; t, x)  \le J_\infty(t, x,  \eta,  Z,  \zeta) - J_n(t, x, \tilde \eta, \tilde Z, \tilde \zeta)\\
 &&\dis =\dbE^{\zeta,Z}\Big[\int_{\t_n^\d}^T [I^h(\zeta_s; s,Z_s) - \eta_s +\eta^{\zeta_{\t^\d_n}; \t^\d_n, R^{\zeta_{\t^\d_n}}_{\t^\d_n}}_s]ds \\
 &&\dis \qq - \sum_{i=1}^n c^P_{\t_i}\1_{\{T_\d< \t_i <T\}} - \sum_{i=n+1}^\infty c^P_{\t_i}\1_{\{\t_i <T\}} -  \sup_\th V^P_0(\th; \t_n)\1_{\{\t_n<T\}}\big]\\
 &&\dis \le \dbE^{\zeta,Z}\Big[\int_{\t_n^\d}^T [C|Z_s| +C_0 - \eta_s]ds + C_0\1_{\{\t_n<T\}}\Big]\\
 &&\dis \le C\dbE^{\zeta,Z}\Big[ \sqrt{(T-\t^\d_n)(\ol L^{\zeta_{\t_n^\d}}_{\t_n^\d} - \cX_{\t_n^\d})}  + (\ol L^{\zeta_{\t_n^\d}}_{\t_n^\d} - \cX_{\t_n^\d}) + \1_{\{\t_n<T\}}\Big]\\
 &&\dis \le C\dbP^{\zeta,Z}(\t_n<T) + C\dbE^{\zeta,Z}\Big[\1_{\{\t_n=T\}} \sqrt{\d(\ol L^{\zeta_{T_\d}}_{T_\d} - \cX_{T_\d})}  + (\ol L^{\zeta_{T_\d}}_{T_\d} - \cX_{T_\d})\Big]
 \eeaa
 where the thrid inequality thanks to \reff{Zetaest1}. Send $\d\to 0$ and note that $\cX_{T-} = 0$, we have
\beaa
J_\infty(t, x,  \eta,  Z,  \zeta) - u_n(\th; t, x) \le C\dbP^{\zeta,Z}(\t_n<T)  \le {C\over n}.
 \eeaa
This implies $u_\infty(\th; t, x) - u_n(\th; t,z) \le  {C\over n}$, and hence $|u_n(\th; t, x) - u_\infty(\th; t,z)| \le  {C\over n}$.
\qed

\bs
\no{\bf Proof of Proposition \ref{prop-uSnmon}.}  We first show that
\bea
\label{u1<u0}
u^S_1(\th_0; t, x) \le u^S_0(\th_0; t, x).
\eea
This, together with \reff{uSn} and Proposition \ref{prop-convS}, implies recursively $u^S_n \downarrow u_\infty$ and  $V^{S,P}_n\downarrow V^P_\infty$.

Indeed, in this case by \reff{uSn} we have
\beaa
u^S_1(\th_0; t, x) = \!\!\! \sup_{(\eta, Z)\in \cA^{\th_0}_{t,x}}\!\!\! \dbE^{\th_0,Z}\Big[\int_t^{\t_1} [I^h(\th_0; s,Z_s) - \eta_s] ds + [u^S_0(\th_0;\t_1, R^{\th_0}_{\t_1}) - c^P_{\t_1}]\1_{\{\t_1<T\}} \Big].
\eeaa
On the other hand, similarly to  \reff{u1DPP}  $u^S_0$ satisfies the dynamic programming principle: 
\beaa
u^S_0(\th_0; t, x) = \sup_{(\eta, Z)\in \cA^{\th_0}_{t,x}} \dbE^{\th_0,Z}\Big[\int_t^{\t_1} [I^h(\th_0; s,Z) - \eta_s] ds + u^S_0(\th_0; \t_1, \ul L^{\th_0}_{\t_1}) \1_{\{\t_1<T\}} \Big].
\eeaa
Note that $R^{\th_0}_{\t_1} > R^{\th_0}_{\t_1} - c^{\th_0}_{\t_1} = \ul L^{\th_0}_{\t_1}$ and $u_0$ is decreasing in $x$, then 
\bea
\label{u0quit}
u^S_0(\th_0; \t_1, R^{\th_0}_{\t_1}) - c^P_{\t_1} \le u^S_0(\th_0; \t_1, \ul L^{\th_0}_{\t_1}),
\eea
which implies \reff{u1<u0} immediately, where the loss is due to the costs $c^{\th_0}_{\t_1}$ and $c^P_{\t_1}$.  
\qed

\end{document}